\documentclass[12pt, reqno]{amsart}
\usepackage{amsthm}
\usepackage{amssymb}
\usepackage{amsmath}

\def\1{\hbox{1\kern-.35em\hbox{1}}}

\newtheorem{theorem}{Theorem}[section]
\newtheorem*{theorem*}{Theorem}
\newtheorem{lemma}[theorem]{Lemma}

\newtheorem*{proposition*}{Proposition}

\newtheorem{definition}[theorem]{Definition}

\newtheorem{remark}[theorem]{Remark}

\newtheorem{example}[theorem]{Example}
\newtheorem{convention}[theorem]{Convention}

\numberwithin{equation}{section}

\newcommand{\bea}{\begin{eqnarray}}
\newcommand{\eea}{\end{eqnarray}}
\newcommand{\be}{\begin{eqnarray*}}
\newcommand{\ee}{\end{eqnarray*}}

\newcommand{\Z}{{\mathbb Z}}

\newcommand{\C}{{\mathbb C}}
\def\VB{\vskip4pt}
\def\VE{\vskip2pt}
\def\VV{{\mathcal V}}
\def\II#1#2{[#1,#2]}
\newcommand{\fb}{{\mathfrak b}}
\newcommand{\fg}{{\mathfrak g}}
\newcommand{\fh}{{\mathfrak h}}

\newcommand{\cC}{{\mathcal C}}

\def\SSP{\vs{-4pt}}

\def\kac#1{V_{#1}}
\def\irr#1{L_{#1}}

\def\low#1{\check #1}

\def\a{\alpha}

\def\D{\Delta}
\def\g{\gamma}
\def\G{\Gamma}
\def\l{\lambda}
\def\L{\Lambda}

\def\sc{\scriptstyle}
\def\ssc{\scriptscriptstyle}
\def\dis{\displaystyle}
\def\cl{\centerline}

\def\ol{\overline}
\def\ul{\underline}
\def\wt{\widetilde}
\def\wh{\widehat}

\def\Rar{\Longrightarrow}

\def\D{\Delta}

\def\Lra{\Longleftrightarrow}

\def\bs{\backslash}
\def\hs{\hspace*}

\def\vs{\vspace*}

\def\ni{\noindent}

\def\N{\mathbb{N}}
\def\Z{\mathbb{Z}}

\def\C{\mathbb{C}}

\def\gl{{\mathfrak {gl}}}

\def\es{\epsilon}

\def\cp#1{{\bar #1}}
\def\tau{\eta}
\def\nn{\mbox{\textit{\textbf n}}}
\def\qq{\mbox{\textit{\textbf q}}}
\def\cc{\mbox{\textit{\textbf c}}}
\def\equa#1#2{
\begin{equation}\label{#1}#2\end{equation}}
\def\equan#1#2{$$#2$$}
\numberwithin{equation}{section}
\begin{document}
%
\title[Composition factors of Kac-modules for ${\mathfrak{gl}}_{m|n}$]
{Composition factors of Kac-modules for 
the general linear Lie superalgebras ${\mathfrak{gl}}_{m|n}$}
\author[Yucai Su]{Yucai Su}
\address{
School of Mathematics and Statistics,
University of Sydney, NSW 2006, Australia;
\newline \indent
Department of Mathematics, Shanghai Jiaotong University, Shanghai
200030, China.
} \email{yucai@maths.usyd.edu.au}

\begin{abstract}
The composition factors of Kac-modules for 
the general linear Lie superalgebras ${\mathfrak{gl}}_{m|n}$
are explicitly determined. In particular, a conjecture of Hughes, King and
van der Jeugt in [{\it J.~Math.~Phys.}, {\bf41} (2000), 5064-5087] is proved.
\end{abstract}
\maketitle

%
\section{Introduction}
%
%
Following the classification of simple Lie superalgebras \cite{Kac0, Kac}, Kac studied finite-dimensional
modules of the classical Lie superalgebras \cite{Kac1, Kac2},
distinguishing between typical and atypical
modules. He also introduced what is now called the Kac-module $\kac{\l}$, 
which was shown to be simple if and only if $\l$ is typical. 
Since then,
Kac-modules,
which themselves
encapsulate rich
information on the structure of the representations,
have been playing extremely active roles in the representation theory of
Lie superalgebras.
For $\l$ atypical, 
the structure of $\kac{\l}$,  
or more generally the problem of classifying finite-dimensional indecomposable modules
has been the subject of intensive study (see, e.g., the References).
By analyzing structures of Kac-modules, van der Jeugt \cite{V} constructed a character formula
for all finite-dimensional irreducible modules over the orthosymplectic Lie superalgebras
${\mathfrak{osp}}_{2|2n}$.
However in the case of the general linear
Lie superalgebras $\gl_{m|n}$, it turned out that
the analysis of structures of Kac-modules is a technical and difficult problem.

There were many partial results on describing structures of Kac-modules or determining
character formulae of irreducible modules over $\gl_{m|n}$. However the full problem 
remained open until
Serganova \cite{Se96, Se98}, based on ideas from
Kazhdan-Lusztig theory,
derived an algorithm for computing 
character formulae for irreducible modules $\irr{\l}$, and
determining the multiplicities $a_{\l,\mu}=[\kac{\l}:\irr{\mu}]$ of
composition factors $\irr{\mu}$ of Kac-modules 
$\kac{\l}$.
(The implementation of this algorithm
turned out to be rather unwieldy to use, e.g., a fact
which was conjectured  by van der Jeugt and Zhang \cite{VZ}
and proved by Brundan \cite{B}
that the composition multiplicities $a_{\l,\mu}$ of the Kac-modules are all either $0$ or $1$,
does not seem 
to follow easily from Serganova's formula since that involves certain alternating sums.)
This work was further developed in \cite{B}, where Brundan used quantum group
techniques to develope
a very practicable algorithm for computing Kazhdan-Lusztig polynomials  
for finite-dimensional irreducible modules over $\gl_{m|n}$
and proved
a theorem previously conjectured in \cite{VZ},
which determines all weights $\l$ such that $a_{\l,\mu}=1$, for a given $\mu$
(there are precisely $2^r$ such $\l$'s, where $r$ is the 
degree of atypicality of $\mu$).
This algorithm was further implemented by Zhang and the author \cite{SZ1}, who 
obtained some closed formulae to compute Kazhdan-Lusztig polynomials, characters and
dimensions for all finite-dimensional irreducible modules over $\gl_{m|n}$.

Brundan's result is useful in understanding structures of Kac-modules.
However this result is not ready to be used in describing 
the structure of a given Kac-module as clear as one would wish,
it still seems to be a problem on how to explicitly 
determine the composition factors $\irr{\mu}$
of the Kac-module $\kac{\l}$, for a given $\l$.
Due to the crucial role that Kac-modules have been playing in the representation
theory of Lie superalgebras, it seems to us that it is highly desirable to
derive a closed formula for computing the composition factors of $\kac{\l}$.
Hughes, King and van der Jeugt \cite{HKV} described an algorithm to determine all
the composition factors of Kac-modules for ${\mathfrak{gl}}_{m|n}$. 
They conjectured that 
there exists a bijection between the composition
factors of $\kac{\l}$ and certain permissible codes (see Definition \ref{code}). 
This conjecture,
which withstood extensive tests against computer
calculations for a wide range of weights $\l$,
describes clearly the structure of $\kac{\l}$.

In this paper, we shall further implement Brundan's result to determine
explicitly the composition factors of the finite-dimensional Kac-modules over the general
linear Lie superalgebras.
A closed formula is obtained for determining the set of the composition factors 
of $\kac{\l}$ (see Theorem
\ref{theo-main}).
This result is quite explicit and easy to apply. In particular we are able to prove
the conjecture of Hughes et al (see Theorem \ref{code-primitive}).
The techniques used in the paper are  
purely combinatorial.

The organization of the paper is as follows. Some background
material on $\gl_{m|n}$ which will be used in the paper is recalled in Section 2. 
In Section 3,
the notion of $\nn\qq\cc$-relationship is introduced that is crucial 
in the proof of the main theorem, which is also
presented in this section. Section 4 is devoted to a proof of the conjecture
of Hughes et al after the notion of permissible codes being introduced, 
and the final section is devoted to the proof of the main result.
Finally we may like to mention that, as is stated earlier,
due to the fact that a Kac-module itself has a complicated
structure, some arguments in the proof may render technical.
\section{Preliminaries}
\label{Preliminaries}
%
%
%
Denote by $\C^{m|n}$ the $\Z_2$-graded vector space with even subspace
$\C^m$ and odd subspace $\C^n$. Then ${\rm End}_{\C}(\C^{m|n})$ with
the $\Z_2$-graded commutator forms the {\it general linear superalgebra}
$\gl_{m|n}$,
which is denoted by $\fg$ throughout the paper.
Choose a basis
$\{ v_a \,\,|\,\, a\in {\bf I}\}$,
for $\C^{m|n}$, where ${\bf I}=\{1, 2, \ldots , m+n\}$, and $v_a$ is
even if $a\le m$, and odd otherwise. 
Let $E_{a b}$ be the matrix unit, namely, the $(m+n)\times(m+n)$-matrix
with all entries being zero except that at the $(a, b)$ position which is $1$.
Then $\{E_{a b}\, |\,a,b\in{\bf I}\}$ forms a basis of $\fg$,
with $E_{a b}$ being even if $a, b\le m$, or
$a, b> m$, and odd otherwise. 
Define the map 
\equan{map}
{[\ ]: {\bf I}\rightarrow \Z_2,\ \
[a]=\big\{\begin{array}{l l}
               \bar{0},  & \mbox{if} \ a\le m, \\
               \bar{1},  & \mbox{if} \ a>m.
              \end{array} 
}
Then the commutation relations 
can be written as
\be [E_{a b}, \ E_{c d}] &=& E_{a d}\delta_{b c}
 - (-1)^{([a]-[b])([c]-[d])} E_{c b}\delta_{a d}.
\ee

The upper triangular matrices form a {\it Borel subalgebra} $\fb$ of
$\fg$,
which contains the {\it Cartan subalgebra} $\fh$ of diagonal matrices.
Let $\{\epsilon_a \,\,|\,\, a\in{\bf I}\}$ be the basis
of $\fh^*$ such that $\epsilon_a(E_{b b})=\delta_{a b}$.
The supertrace induces a {\em bilinear form} $(\;
, \: ): \fh^*\times \fh^* \rightarrow \C$ on $\fh^*$ such
that $(\epsilon_a, \epsilon_b)=(-1)^{[a]} \delta_{a b}$.
Relative to the Borel subalgebra $\fb$,
the roots of $\fg$
can be expressed as $\epsilon_a-\epsilon_b, \:\, a\ne b$, where
$\epsilon_a-\epsilon_b$ is even if $[a]+[b] = \bar{0}$ and odd
otherwise. The set of the {\em positive roots}  is $\D^+=\{
\epsilon_a-\epsilon_b \,\,|\,\, a< b\}$, and the set of {\em simple roots}
is $\{\epsilon_a-\epsilon_{a+1} \,\,|\,\, a<m+n\}$.

We denote ${\bf I}^1=\{1,2,...,m\}$ and
${\bf I}^2=\{\cp 1,\cp2,...,\cp n\}$, where here and below we use
the notation $$\cp \nu=\nu+m.$$
Then ${\bf I}={\bf I}^1\cup{\bf I}^2$.
The sets of {\em positive
even roots} and {\em odd roots} are respectively
\begin{eqnarray*}
\D_0^+&\!\!\!=\!\!\!&\{\a_{i,j}=\es_i-\es_j,\,
\a_{\nu,\eta}=\es_\nu-\es_\eta
\,\,|\,\,1\le i<j\le m,\,\cp 1\le\nu<\eta\le\cp n\},
\\
\D_1^+&\!\!\!=\!\!\!&\{\a_{i,\nu}=\es_i-\es_\nu\,\,|\,\,i\in{\bf I}^1,\,
\nu\in{\bf I}^2\}.
\end{eqnarray*}
The Lie algebra $\fg$ admits a $\Z_2$-consistent $\Z$-grading
$$\fg=\fg_{-1}\oplus\fg_0\oplus\fg_{+1},
\mbox{ where }\fg_0=\fg_{\bar0}\cong\gl(m)\oplus\gl(n)
\mbox{  and }\fg_{\pm1}\subset\fg_{\bar1},
$$
with
$\fg_{+1}$ (resp.~$\fg_{-1}$) being the nilpotent subalgebra spanned by the odd positive (resp.~negative) root
spaces.
We define a total order on $\D_1^+$ by
\equan{0.0}
{
\a_{i,\nu}<\a_{j,\eta}\,\ \ \Lra\,\ \
\nu-i<\eta-j\mbox{ or }\nu-i=\eta-j\mbox{ but }i>j.
}

An element in $\fh^*$ is called a {\em weight}. A weight $\L$ is {\em integral} if
$(\Lambda, \, \es_a)\in\Z$ for all $a$,
 and {\em dominant} if
${2(\Lambda, \, \alpha)}/{(\alpha,\, \alpha)}\ge0$ for all positive even roots $\a$ of $\fg$.
Denote by $P$ (resp.~$P_+$) the set of integral (resp.~dominant integral)
weights.

Since $\{\es_i\,\,|\,\,i\in{\bf I}\}$ is a $\C$-basis of $\fh^*$, a weight $\l\in\fh^*$ can be
written as $\l=\sum_{i\in{\bf I}}\l'_i\es_i$ with $\l'_i\in\C$, 
and it is usually denoted by
\equa{weight0}
{\mbox{$
\l=(\l'_1,\l'_2,...,\l'_m\,\,|\,\,\l'_{\cp 1},\l'_{\cp2},...,\l'_{\cp n}).
$}}
But sometimes, we shall find it is more convenient to denote the
weight $\l$ by
\equa{weight1}
{\mbox{$
\l=(\l_1,\l_2,...,\l_m\,\,||\,\,\l_{\cp 1},\l_{\cp2},...,\l_{\cp n}),
\mbox{ where }\l_i=\biggl\{\begin{array}{cl}\l'_i+i&\mbox{if }i\in{\bf I}_1,\\[2pt]
-\l'_i+i-m&\mbox{if }i\in{\bf I}_2.\end{array}
$}}
One can easily convert notation (\ref{weight0}) to notation (\ref{weight1}), or vice versa.
With notation (\ref{weight1}), 
the set of integral weights coincides with
the set $\Z^{m|n}$ of $(m+n)$-tuples of integers
and the set of dominant integral weights coincides with the subset
$\Z_+^{m|n}$  of $(m+n)$-tuples $\l$ satisfying 
\equa{d-weight}
{
\l_1>\l_2>\cdots>\l_m,\ \ \l_{\cp1}<\l_{\cp2}<\cdots<\l_{\cp n}.
}
We should note that there is no loss of generality in restricting our attention to
integral weights $\l$ since an arbitrary ``$r$-fold atypical''
(see below)
finite-dimensional Kac-module can be obtained from $\kac{\l}$, for some $\l\in\Z^{m|n}$,
by tensoring with one-dimensional module.

Let $W=S_m\times S_n$ be the {\em Weyl group} of $\fg$, where $S_m$ is the
{\it symmetric group}
of degree $m$. The action of $W$ on $P$, by definition, is 
\equa{action-of-W}
{w\l=(\l_{w(1)},\l_{w(2)},...,\l_{w(m)}\,\,||\,\,
\l_{w(\cp1)},\l_{w(\cp2)},...,\l_{w(\cp n)})\in\Z^{m|n}\mbox{ for }\l\in\Z^{m|n},
}
where $w\in W$.
An integral weight $\l$ is called {\it regular}
or {\it non-vanishing} (in sense of \cite{HKV, VHKT}) if it is $W$-conjugate 
to a dominant weight
(which is denoted by $\l^+$ throughout the paper),
otherwise it is called {\it vanishing}.

For a regular weight $\l$ in (\ref{weight1}),
we define the {\it atypicality matrix of $\l$}
to be the $m\times n$ matrix  
\equa{Aty-matrix}
{
A(\l)=(A(\l)_{i,\eta})_{m\times n}, \mbox{ where }
A(\l)_{i,\eta}=\l_i-\l_{\cp \eta},\ 1\le i\le m,\
1\le \eta\le n.
}
An odd root $\a_{i,\cp\eta}$ is an {\it atypical root of $\l$} if
$A(\l)_{i,\eta}=0$. Let $\G_\l=\{\a_{i,\cp\eta}\,\,|\,\,A(\l)_{i,\eta}=0\}$
be the set of atypical roots, and
$r=\#\G_\l$ be the {\it degree of atypicality}. 
We also denote $\#\l=r$. 
Thus $\#\l$ is the number of pairs $(i,\cp\eta)$ whose entries are equal: $\l_i=\l_{\cp\eta}$.
A weight $\l$
is {\it typical} if $r=0$; {\it atypical} if $r>0$ (in this case $\l$ is also called
an {\em $r$-fold atypical weight}). 
If $\l$ is dominant and $r$-fold atypical,
we label its atypical roots by $\g_1,...,\g_r$ ordered in such a way that (cf.~Example \ref{exam3.1})
\equa{aty-roots}
{
\biggl\{\begin{array}{l}
\g_1<\g_2<\cdots<\g_r
\mbox{ and }\g_s=\a_{m_s,n_s},\ \ s=1,2,...,r,\mbox{ with}\\[4pt]
1\le m_r<m_{r-1}<\cdots<m_1\le m<n_1<n_2<\cdots<n_r\le\cp n.
\end{array}
}

For an integral dominant weight $\l$, denote by $\irr{\l}^{(0)}$ the
finite-dimensional
irreducible $\fg_0$-module with highest weight $\l$. Extend it to a
$\fg_0\oplus\fg_{+1}$-module by putting $\fg_{+1}\irr{\l}^{(0)}=0$. Then
the {\em Kac-module} $\kac{\l}$ is the induced
module
\equan{Kac-module}
{
\kac{\l}={\rm Ind}_{\fg_0\oplus\fg_{+1}}^{\fg}\irr{\l}^{(0)}\cong U(\fg_{-1})
\otimes_{\C}\irr{\l}^{(0)}.
}
Denote by $\irr{\l}$ the {\em irreducible
module} with highest weight $\l$ (which is the unique irreducible quotient
module of $\kac{\l}$).

The following result is due to Kac \cite{Kac1, Kac2}.
\begin{theorem}
The finite-dimensional Kac-module $\kac{\l}$ over $\fg$ is irreducible if and only if
$\l$ is typical.
\end{theorem}
\section{Composition factors of Kac-modules}
\label{composition}
\def\ol{\bar}
For convenience, we introduce the notation
\equan{I-st}
{
\II{s}{t}=\biggl\{\begin{array}{cl}
\{i\in\Z\,\,|\,\,s\le i\le t\}&\mbox{if }s\le t,
\\[4pt]
\emptyset&\mbox{otherwise},\end{array}
}
for $s,t\in\Z$ (this notation will not be confused 
with the Lie bracket since we do not need to
use that below).
\subsection{The $\nn\qq\cc$-relationship}
Let $P_r$ be the subset of $\Z^{m|n}$ consisting of regular $r$-fold atypical
weights $\l$ such that the atypical roots of $\l$ satisfy (\ref{aty-roots}) and
\equa{i-j}
{\l_i>\l_j, \  i<j\le m, \ i,j\in{\bf I}\bs{\bf I}^\l_r\ \mbox{ and }\
\l_a<\l_b, \ m<a<b,\ a,b\in{\bf I}\bs{\bf I}^\l_r,
}
where ${\bf I}^\l_r=\{m_r,m_{r-1},...,m_1,n_1,n_2,...,n_r\}$.
We denote by $D_r$ the subset of $P_r$ of the elements
$\l$ such that $\l_{m_1}<\l_{m_2}<...<\l_{m_r}$.

For $\l\in P_r$, we denote
\equa{denote-t-f}
{t^\l\in\Z^{r|r}\mbox{ \ or \ }\ol t{\ssc\,}^\l\in\Z^{m-r|n-r}
}
to be the element obtained
from $\l$ by deleting its $i$-th entry for $i\in{\bf I}\bs{\bf I}^\l_r$ or $i\in{\bf I}^\l_r$
respectively.
Thus $\ol t{\ssc\,}^\l$ is always dominant for all $\l\in P_r$, 
and $t^\l$ is dominant if and only if $\l\in D_r$.
We also introduce the following three sets of integers:
\equa{set-R(f)}{
S(\l)\!=\!\{\l_i\,\,|\,\,i\in {\bf I}\},\ 
T(\l)\!=\!S(t^\l)\!=\!\{\l_{m_s}\,\,|\,\,s\in\II{1}{r}\},\ 
\ol T(\l)\!=\!S(\ol t{\ssc\,}^\l)\!=\!S(\l)\bs T(\l).
}
\begin{convention}\label{conv} 
We usually use the superscript `$\l$' to indicate that
a notation is associated with $\l$, like in ${\bf I}^\l_r,\,t^\l,\,\ol t{\ssc\,}^\l$.
However, when confusion is unlikely
to occur, the superscript will be dropped, like in $\ell_{s,t},\,c_{s,t},\,k_i,\,\low k_i$ 
below.
\end{convention}
The following notion of $\nn\qq\cc$-relationship was first introduced by Hughes et al \cite{HKV}
from a different point of view.
\begin{definition}\label{defi3.1} 
Suppose $\l\in D_r$. 
Let $\nn,\qq,\cc$ be three symbols.
For $1\le s\le t\le r$, we define
\equa{c'-st}
{
\ell_{s,t}=\#(\II{f_{m_s}}{f_{m_t}}\bs S(f)),\ \ \
c_{s,t}=\left\{\begin{array}{ll}
\nn&\mbox{if \ \ }\ell_{s,t}>t-s,\\[4pt]
\qq&\mbox{if \ \ }\ell_{s,t}=t-s,\\[4pt]
\cc&\mbox{if \ \ }\ell_{s,t}<t-s.
\end{array}\right.
}
Note that in $(\ref{c'-st})$, we abused the notation by using
$\II{f_{m_s}}{f_{m_t}}\bs S(f)$ to denote
$$\II{f_{m_s}}{f_{m_t}}\bs(\II{f_{m_s}}{f_{m_t}}\cap S(f)).$$
Two atypical roots $\g_s,\g_t$ of $\l$ are called $($cf.~\cite{HKV, VHKT}$)$
\equan{nooo}{
\begin{array}{rlll}
\mbox{\rm(i)}\!\!\!\!&\mbox{ {\it normally related} or {\it $\nn$-related}}&\Lra&
c_{s,t}=\nn;\\[4pt]
\mbox{\rm(ii)}\!\!\!\!&\mbox{ {\it quasi-critically related} or {\it $\qq$-related}}&\Lra&
c_{s,t}=\qq;\\[4pt]
\mbox{\rm(iii)}\!\!\!\!&\mbox{ {\it critically related} or {\it $\cc$-related}}&\Lra&
c_{s,t}=\cc.
\end{array}
}
\end{definition}
Note that $\qq$-relationship is reflexive
and transitive but not symmetric ($c_{t,s}$ is not
defined when $t>s$); $\cc$-relationship or $\nn$-relationship is transitive.
\begin{example}\label{exam3.1}
Suppose
\equa{exam1}
{
\l=\bigl(\stackrel{^{\sc4}}{{15}},11,\stackrel{^{\sc3}}{{10}},\stackrel{^{\sc2}}{7},6,4,
\stackrel{^{\sc1}}{3}\,\bigl|\!\bigr|\,\stackrel{^{\sc1}}{3},5,
\stackrel{^{\sc2}}{7},8,\stackrel{^{\sc3}}{{10}},\stackrel{^{\sc4}}{15}\bigr),\ \,
S(\l)=\{\stackrel{^{\sc1}}{3},4,5,6,
\stackrel{^{\sc2}}{7},8,\stackrel{^{\sc3}}{{10}},11,\stackrel{^{\sc4}}{{15}}\},
}
where we put a label $s$ over an entry to indicate it corresponds to
the $s$-th atypical root $($such an entry is called an atypical entry$)$. Then
\par\ni\hs{8ex}
$c_{1,2}=\cc$
since $\II{3}{7}\bs S(\l)=\emptyset$ is of cardinality $0<1$;
\\\hs{8ex}
$c_{1,3}=\cc$ since  $\II{3}{10}\bs S(\l)=\{9\}$ is of cardinality
$1<2$;
\\\hs{8ex}
$c_{1,4}=\nn$ since $\II{3}{15}\bs S(\l)=\{9,12,13,14\}$ is of cardinality $4>3$;
and
\\\hs{8ex}
$c_{2,3}=\qq,\,c_{2,4}=\nn,\,c_{3,4}=\nn$.
\end{example}
A simple way to determine $c_{s,t}$ is to count the number of integers
between $\l_{m_s}$ and $\l_{m_t}$ which do not belong to the set $S(\l)$.
If the number (which is $\ell_{s,t}$) is smaller than (resp.~equal to,
or bigger than) $t-s$ then $c_{s,t}=\cc$ (resp.~$\qq$, or $\nn$).

The following discussion may illustrate the significance of the concept of
$\nn\qq\cc$-rela\-tion\-ship.

\subsection{Raising and lowering operators}
Let $\l\in D_r$. For $s=1,2,...,r$, we set (recall Convention \ref{conv})
\begin{eqnarray}
\label{p-s}
p_s&\!\!\!=\!\!\!&
\biggl\{\begin{array}{l}s\mbox{ \ if \ }s=r\mbox{ \ or \ } c_{s,s+1}\ne \cc,\\[2pt]
\max\{p\in\II{s+1}{r}\,\,|\,\,\,
c_{s,s+1}=c_{s,s+2}=...=c_{s,p}=\cc\}\mbox{ \ otherwise},\end{array}
\\
\label{p--s}
\low p_s&\!\!\!=\!\!\!&
\biggl\{\begin{array}{l}s\mbox{ \ if \ }s=1\mbox{ \ or \ } c_{s-1,s}\ne \cc,\\[2pt]
\min\{p\in\II{1}{s-1}\,\,|\,\,\,c_{p,s}=c_{p+1,s}=...=c_{s-1,s}=\cc\}
\mbox{ \ otherwise}.\end{array}
\end{eqnarray}
Namely, $p_s$ (resp.~$\low p_s$) is the largest (resp.~smallest)
integer such that all $\g_i$ with $i\in\II{s+1}{p_s}$ (resp.~$i\in\II{\low p_s}{s-1}$)
are $\cc$-related to $\g_s$.
Define $r$-tuples $(k_1,k_2,...,k_r)$ and 
$(\low k^{(\nu)}_1,\low k^{(\nu)}_2,...,\low k^{(\nu)}_r)$, $\nu\ge1$, 
of positive integers associated with $\l\in D_r$ by:
\begin{eqnarray}
\label{k-i}
k_s&\!\!\!=\!\!\!&\min
\left\{k>0\,\left|\ \#\left({}^{^{^{^{\,\!}}}}\II{\l_{m_s}}{\l_{m_s}+k}\bs S(\l)\right)=p_s+1-s\right\},\right.
\\
\label{k--i}
\low k^{(\nu)}_s&\!\!\!=\!\!\!&\min
\left\{k>0\,\left|\ \#
\left({}^{^{^{^{\,\!}}}}
\II{\l_{m_s}-k}{\l_{m_s}}\bs S(\l)\right)=\nu\right\}.\right.
\end{eqnarray}
This definition means that $\l_{m_s}+k_s$ (resp.~$\l_{m_s}-\low k^{(\nu)}_s$) is the
$(p_s+1-s)$-th smallest (resp.~the $\nu$-th largest) integer not in the set $S(\l)$ which is bigger (resp.~smaller)
than $\l_{m_s}$.
For convenience we set $\low k_i^{(0)}=0$ and
we simply denote $\low k_i=\low k^{(1)}_i$ for $i=1,2,...,r$.

A simple general way to compute $k_i$ is the following procedure:
First set $S=S(\l)$. Suppose we have computed $k_r,k_{r-1},...,k_{i+1}$.
To compute $k_i$, we count the integers in the set $S$ starting with $\l_{m_i}$
and stop at the first integer, say $k$, not in $S$. 
Then $k_i=k-\l_{m_i}$. {\bf Now add $k$ into the set $S$}, and continue.

The computation of $\low k^{(\nu)}_i$ is much simpler:
count the integers downward in the set $S(\l)$ 
starting with $\l_{m_i}$ until we find $\nu$ integers not in $S(\l)$. 
Say we stop at the integer $\low k$, then $\low k^{(\nu)}_i=\l_{m_i}-\low k$.

For example, if $\l$ is as in $(\ref{exam1})$, then
\equa{exam-k-i}{
(k_1,k_2,k_3,k_4)=(10,2,2,1)\mbox{ \ and \ }(\low k_1,\low k_2,
\low k_3,\low k_4)=(1,5,1,1).
}

If $\l\in P_r$ (not necessarily in $D_r$),
we can still compute $k_i$ and $\low k^{(\nu)}_i$ in the above way, but the difference
lies in that the
$k_i$'s are computed in the order that
each time we compute $k_i$ with $\l_{m_i}$ being the largest 
among all those
$\l_{m_i}$'s, the corresponding $k_i$'s of which are not yet computed.
\begin{lemma}\label{lemm3.1}
Suppose $\l\in D_r$.
\par
{\rm(1)}
For each $s\in\{1,2,...,r\}$, $k_s$
is the smallest positive integer such that
\begin{eqnarray}
\label{cond3.1}
\!\!\!\!&\!\!\!\!&\!\!\!\!
(\l+\theta_t k_td_{m_t,n_t})+k_sd_{m_s,n_s}\mbox{ is regular for }
\theta_t\in\{0,1\},\, s<t\le r,
\end{eqnarray}
where $d_{i,j}\in\Z^{m|n}$ is the $(m+n)$-tuple whose entries are zero except the
$i$-th and $j$-th entries which are $1$.

{\rm(2)}
For each $s\in\{1,2,...,r\}$,
$\low k_s$ is the smallest positive integer such that
\equa{k'}{\mbox{
$\l-\low k_sd_{m_s,n_s}\mbox{ is regular}.$ 
}}

{\rm(3)}
The tuple $(k_r,k_{r-1},...,k_1)$
is the lexicographically
smallest tuple of positive integers such that for all $\theta=(\theta_1,\theta_2,...,
\theta_r)\in\{0,1\}^r$, $\l+\sum_{s=1}^r\theta_sk_sd_{m_s,n_s}$ 
is regular. Thus
$(k_r,k_{r-1},...,k_1)$ is the tuple satisfying \cite[{\rm Main Theorem}]{B}.
\end{lemma}
\begin{proof} See \cite[Lemma 3.3]{SZ1}.
It can also be proved directly using the definitions (\ref{k-i}) and (\ref{k--i}).
\end{proof}
Following \cite{B} (see also \cite{VZ, SZ1}), we define the {\it raising operator}
$R_{m_s,n_s}$ and the {\it lowering operator} $L_{m_s,n_s}$ on $P_r$ by
\equa{raising}
{
R_{m_s,n_s}(\l)=\l+k_sd_{m_s,n_s},
\ \ 
L_{m_s,n_s}(\l)=\l-\low k_sd_{m_s,n_s},
}
for $s=1,2,...,r$ and $\l\in P_r$.
It is straightforward to verify that the composition of $\nu$ copies of 
the operator $L_{m_s,n_s}$ is
\equa{low-nu}{
L^\nu_{m_s,n_s}(\l)=\l-\low k^{(\nu)}_sd_{m_s,n_s}\mbox{ for }\nu\ge1.
}
Let $\theta=(\theta_1,\theta_2,...,\theta_r)\in\N^r$, where $\N=\{0,1,...\}$.
We define
\begin{eqnarray}
\label{R'-theta-f}
\!\!\!\!&\!\!\!\!&\!\!\!\!
R_\theta(\l)=(R^{\theta_r}_{m_r,n_r}\circ R^{\theta_{r-1}}_{m_{r-1},n_{r-1}}\circ\cdots\circ R^{\theta_1}_{m_1,n_1}(\l))^+,
\\
\label{l-theta-f}
\!\!\!\!&\!\!\!\!&\!\!\!\!
L'_\theta(\l){\ssc\,}=(L^{\theta_r}_{m_r,n_r}{\ssc\,}\circ L^{\theta_{r-1}}_{m_{r-1},n_{r-1}}{\ssc\,}\circ\cdots\circ L^{\theta_1}_{m_1,n_1}{\ssc\,}(\l))^+,
\end{eqnarray}
where in general $\l^+$ denotes the unique dominant element which is $W$-conjugate to $\l$.

\begin{remark}\label{rema3.0}
We remark
that the definition of lowering operator defined here is different from that defined
in \cite[\S3-f]{B}. In \cite[\S3-f]{B},  $\low k_s$ appearing in $(\ref{raising})$
needs to satisfy the following condition instead of condition $(\ref{k'})$:
\equan{con-there}{\mbox{
$(\l-\theta_p\low k_pd_{m_p,n_p})-\low k_sd_{m_s,n_s}$ is regular for all $1\le p<s,\,
\theta_p\in\{0,1\}$.
}
}
\end{remark}
From definitions (\ref{k-i}) and (\ref{R'-theta-f}), by induction
on $\#\{s\in\II{1}{r}\,\,|\,\,\theta_s=1\}$, we can prove (cf.~\cite[Main Theorem]{B}
or Lemma \ref{lemm3.1}(3))
\equa{Main-t}
{
R_\theta(\l)=\bigl(\l+\mbox{$\sum\limits_{s=1}^r$}\theta_sk_sd_{m_s,n_s}\bigr)^+
\mbox{ for }\theta\in\{0,1\}^r.
}
Also from (\ref{l-theta-f}) and (\ref{low-nu}), by induction
on $\#\{s\in\II{1}{r}\,\,|\,\,\theta_s\ne0\}$, we can prove
\equa{Main--t}
{
L'_\theta(\l)=\bigl(\l-\mbox{$\sum\limits_{s=1}^r$}\low k\!\!\!k\!\!\!k^{(\theta_s)}_sd_{m_s,n_s}\bigr)^+
\mbox{ for }\theta\in\N^r,
}
where $\low k\!\!\!k\!\!\!k^{(\theta_s)}_s$ is define as follows: first set
\equa{set-l-s-}
{\l^{(0)}=\l\mbox{ and }
\l^{(s)}=L^{\theta_s}_{m_s,n_s}\circ\cdots\circ L^{\theta_1}_{m_1,n_1}(\l)\mbox{ if }
s\in\II{1}{r},
}
then $\low k\!\!\!k\!\!\!k^{(\theta_s)}_s$ is $(\low k_s^{(\theta_s)})^{\l^{(s-1)}}$ 
(which is $\low k_s^{(\theta_s)}$ defined by $\l^{(s-1)}$, cf.~Convention \ref{conv}).

One can also observe that
\equa{l-mu}
{
\l=L'_\theta(\mu)\mbox{ or }\l=R_\theta(\mu)\ \ \Rar\ \ 
\ol t{\ssc\,}^\l=\ol t{\ssc\,}^\mu
\mbox{ \ (cf.~(\ref{exam2}))}.
}
\subsection{Composition factors of Kac-modules}
Let $\l,\mu\in\Z_+^{m|n}$. 
We define $a_{\l,\mu}=$ $[\kac{\l},\irr{\mu}]$ to be the multiplicity
of irreducible module $\irr{\mu}$ in the composition series of 
the Kac-module $\kac{\l}$. If $a_{\l,\mu}\ne0$, we say $\irr{\mu}$ is a {\it composition
factor} of $\kac{\l}$, and $\mu$ is a {\it primitive weight} of $\kac{\l}$.

The following theorem, which was conjectured in \cite[Conjecture 4.1]{VZ},
is due to Brundan \cite[Main Theorem]{B}.
\begin{theorem}\label{theo3.1}
Let $\mu\in\Z_+^{m|n}$ be an $r$-fold atypical dominant integral weight. Then
$a_{\l,\mu}\le1$ for all $\l$ and 
\equa{theo1-1}
{a_{\l,\mu}=1\ \ \Lra\ \ 
\l=R_\theta(\mu)\mbox{ for some }\theta\in\{0,1\}^r.}
\end{theorem}
As stated in the introduction, this theorem is useful in understanding structures
of Kac-modules, but it is still desirable to give a closed formula to compute 
the composition factors of the Kac-module $\kac{\l}$, for a given $\l$.
Below we shall implement this theorem to derive such a formula (see Theorem \ref{theo-main}).

Denote 
\equa{Theta-r}
{
\Theta_r=\{\theta=(\theta_1,\theta_2,...,\theta_r)\in\N^r\,\,|\,\,\theta_s\le s
\mbox{ for }s=1,2,...,r\},
}
a subset of $\N^r$ of cardinality $(r+1)!$.
\begin{definition}\label{defi3.4}
Let $\l\in D_r$. We define $\Theta^\l$ to be the subset of $\Theta_r$ consisting of $\theta=(\theta_1,\theta_2,...,\theta_r)$
satisfying the following conditions.

For $s=1,2,...,r$, if $\theta_s\ne0$, then
\equa{conds1}
{
c_{s-\theta_s,s}\ne \cc,\mbox{ and } c_{s+1-\theta_s,s}\ne \nn,
}
and
for all $p\in\II{s+1-\theta_s}{s-1}$,
\equa{conds2}
{\theta_p\le \theta_s-s+p,\mbox{ and equality implies } c_{p,s}=\cc,
}
and furthermore in case $c_{p,s}\ne \nn$,
\equa{conds4}
{
\theta_p\ne0,\mbox{ or }
\exists\,p'\in\II{p+1}{s}\mbox{ such that }c_{p,p'}=\qq\mbox{ and }
\theta_{p'}\ge p'+1-p,
}
\equa{conds5}
{
\mbox{with the exception that }
p=s+1-\theta_s<s\mbox{ and }c_{p,s}=\qq\, \ \Rar\ \,
\theta_p=0.
}
\def\NOTHING{
\begin{eqnarray}
\label{conds2}
\!\!\!\!\!\!&\!\!\!\!\!\!&\!\!\!\!\!\!\!\!\!\!\!\!\!\!\!\!\!\!
\theta_p\le \theta_s-s+p,
\\
\label{conds4}
\!\!\!\!\!\!&\!\!\!\!\!\!&\!\!\!\!\!\!\!\!\!\!\!\!\!\!\!\!\!\!
c_{p,s}\ne n\, \Rar\,
\theta_p\ne0,\mbox{ or }
\exists\,p':\:p{\sc\!}<{\sc\!}p'{\sc\!}\le{\sc\!} s,\,c_{p,p'}{\sc\!}={\sc\!}q
\mbox{ and }\theta_{p'}{\sc\!}\ge{\sc\!} p'{\sc\!}+{\sc\!}1{\sc\!}-{\sc\!}p
\mbox{ with}
\\
\label{conds5}
\!\!\!\!\!\!&\!\!\!\!\!\!&\!\!\!\!\!\!\!\!\!\!\!\!\!\!\!\!\!\!
\phantom{c_{p,s}\ne n\, \Rar\, }
\mbox{the exception that }
p{\sc\!}={\sc\!}s{\sc\!}+{\sc\!}1{\sc\!}-{\sc\!}\theta_s\mbox{ and }
c_{p,s}{\sc\!}={\sc\!}q\mbox{ implies }\theta_p{\sc\!}={\sc\!}0.
\end{eqnarray}
}
\end{definition}
\begin{convention}\label{conv2} 
We use the convention that if
an undefined notation appears in an expression, then this expression
is omitted; for instance in $(\ref{conds1})$
condition $c_{s-\theta_s,s}\ne \cc$ is omitted in case $\theta_s=s$.
\end{convention}
\begin{remark}\label{rema3.1}
{\rm(1)} If $c_{s,t}=\cc$ for all $s<t$ $($in this case the atypical roots of $\l$ are called
totally $\cc$-related$)$,
then the first condition of
$(\ref{conds1})$ implies  $s=\theta_s$ and
condition $(\ref{conds4})$ implies $\theta_p\ne0$. Thus 
\equa{total-c}
{
\Theta^\l=\{(1,2,...,p,0,...,0)\,\,|\,\,p=0,1,...,r\}
\mbox{ is of cardinality }r+1.
}

{\rm(2)} If $c_{s,t}=\nn$ for all $s<t$  $($in this case the atypical roots of $\l$ are called
totally $\nn$-related$)$,
then the second condition of $(\ref{conds1})$ implies 
$s+1-\theta_s=s$, i.e., $\theta_s=1$, and other conditions are all missing since $p$ does not exist in this case. So
\equa{total-n}
{
\Theta^\l=\{0,1\}^r
\mbox{ is of cardinality }2^r.
}
\end{remark}
Let us describe how to determine all elements of $\Theta^\l$ in general.
To do this, we need some more notations.
For $1\le s\le t\le r$, denote by (roughly speaking, $\l^{(s,t)}$ defined below
only keeps those entries of $\l$ ranging from the $s$-th atypical entry to
the $t$-th atypical entry)
\equa{l(s,t)}
{
\l^{(s,t)}
}
the $(t-s+1)$-fold atypical weight
(for some Lie superalgebra ${\mathfrak {gl}}_{k|l}$ with $k\le m,\,l\le n$) obtained
from $\l$ by deleting all entries $\l_i$ for
$1\le i<m_t$ or $m_{s-1}<i\le m$ (set $m_{0}=m+1$), and 
deleting $\l_{\cp\eta}$ for
$m<\cp\eta<n_{s-1}$ (set $n_{0}=m$)
or $n_t<\cp\eta\le\cp n$. 
For instance, if $\l$ is as in
(\ref{exam1}), then 
\equa{l(2,3)}{
\l^{(2,3)}=
(\ul{10},\ul{7},6,4\,||\,5,\ul{7},8,\ul{10}),\ \ \
\l^{(3,4)}=
(\ul{15},11,\ul{10}\,||\,8,\ul{10},\ul{15}),
}
where here and below, the underlined entries
are the atypical entries, i.e., entries
corresponding to the atypical roots.
Then we can define $\Theta^{\l^{(s,t)}}\subset\Theta_{t-s+1}$ by Definition 
\ref{defi3.4}.
Let $\theta\in\N^r$ and $S\subset\II{1}{r}$. We denote
\equa{theta-S}
{
\theta_S\in\N^{\#S}
}
to be the element obtained from $\theta$ by deleting entries $\theta_i$ for $i\notin S$.

By Definition \ref{defi3.4}, we immediately obtain the following lemma which
can be used to determine all elements of $\Theta^\l$ 
(by the procedure of induction on $r$).
\begin{lemma}\label{lemm3.3}
Let 
\equa{M}
{M=\{s\in\II{1}{r+1}\,\,|\,\,c_{s-1,r}\ne \cc,\,c_{s,r}\ne \nn\},
}
$($note that Convention \ref{conv2} means that $r+1\in M)$,
and for $s\in M$ denote
\equa{denote-theta-l-s}
{
\Theta^{\l,s}=
\left\{\begin{array}{ll}
\{\theta\in\Theta_r\,\,|\,\,
\theta_{\II{2}{r}}\in\Theta^{\l^{(2,r)}},\,\theta_1=0,\,\theta_r
=r\}& \mbox{if }s=1,{\ssc\,}c_{1,r}=\qq,
\\[2pt]
\{\theta\in\Theta_r\,\,|\,\,
\theta_{\II{1}{r-1}}\in\Theta^{\l^{(1,r-1)}},\,\theta_1\ne0,\,\theta_r
=r\}& \mbox{if }s=1,{\ssc\,}c_{1,r}=\cc,
\\[2pt]
\{\theta\in\Theta_r\,\,|\,\,
\theta_{\II{1}{s-1}}\in\Theta^{\l^{(1,s-1)}},
\\[2pt]
\phantom{\{\theta\in\Theta_r\,\,|\,\,}
\theta_{\II{s}{r}}\in\Theta^{\l^{(s,r)}},\,\theta_r
=r+1-s\}&\mbox{if }2\le s\le r,
\\[2pt]
\{\theta\in\Theta_r\,\,|\,\,\theta_{\II{1}{r-1}}\in\Theta^{\l^{(1,r-1)}},\,
\theta_r=0\}& \mbox{if }s=r+1.
\end{array}\right.
}
\def\NOUSE{
\equa{denote-theta-l-s}
{\,\,\!\!
\Theta^{\l,s}\!=\!\!
\left\{\!\!\!\begin{array}{ll}
\{\theta\!\in\!\Theta_r\,\,|\,\,
\theta_{\II{2}{r}}\!\in\!\Theta^{\l^{(2,r)}},\,\theta_1\!=\!0,\,\theta_r
\!=\!r\}&\!\!\! \mbox{if }s\!=\!1,{\ssc\,}c_{1,r}\!=\!q,
\\[2pt]
\{\theta\!\in\!\Theta_r\,\,|\,\,
\theta_{\II{1}{r-1}}\!\in\!\Theta^{\l^{(1,r-1)}},\,\theta_1\!\ne\!0,\,\theta_r
\!=\!r\}&\!\!\! \mbox{if }s\!=\!1,{\ssc\,}c_{1,r}\!=\!c,
\\[2pt]
\{\theta\!\in\!\Theta_r\,\,|\,\,
\theta_{\II{1}{s-1}}\!\in\!\Theta^{\l^{(1,s-1)}},\,
\theta_{\II{s}{r}}\!\in\!\Theta^{\l^{(s,r)}},\,\theta_r
\!=\!r\!+\!1\!-\!s\}&\!\!\!\mbox{if }2\le s\le r,
\\[2pt]
\{\theta\!\in\!\Theta_r\,\,|\,\,\theta_{\II{1}{r-1}}\!\in\!\Theta^{\l^{(1,r-1)}},\,
\theta_r\!=\!0\}&\!\!\! \mbox{if }s=r+1.
\end{array}\right.\!\!\!\!\!\!\!\!\!\!
}
}
Then
\equa{Theta-l-for}
{
\Theta^\l=
\bigcup_{s\in M}\Theta^{\l,s}\mbox{ \ \ $($disjoint union$)$}.
}
\end{lemma}
\begin{example}\label{exam3.2}
Suppose $\l$ is as in $(\ref{exam1})$. Then $\Theta^\l$ contains the following $14$
elements:
\equa{exam2-1}
{\begin{array}{l}
(0,0,0,0),\,(1,0,0,0),\,(1,2,0,0),\,(0,0,1,0),\,(1,0,1,0),\,(1,2,1,0),\,(1,0,3,0),\\[4pt]
(0,0,0,1),\,(1,0,0,1),\,(1,2,0,1),\,(0,0,1,1),\,(1,0,1,1),\,(1,2,1,1),\,(1,0,3,1).
\end{array}
}
\end{example}
The following interesting fact was observed in \cite{HKV}.
\begin{lemma}\label{lemm3.2}
If $c_{s,t}=\qq$ for all $s<t$
$($in this case the atypical roots of $\l$ are called
totally $\qq$-related$)$,
 then
\equa{rema1}
{
\#\Theta^\l=\frac{1}{r+2}
\bigl(\!\begin{array}{c}2r+2\\r+1\end{array}\!\bigr)=C_{r+1},
}
where $C_r=\frac{1}{r+1}({\ssc\,}^{2r}_{\,r}{\ssc\,})$ is the well-known
$r$-th Catalan number.
\end{lemma}
\begin{proof}
It is easy to check (\ref{rema1}) if $r=1,2$. Suppose $r\ge3$. 
Note from (\ref{denote-theta-l-s}) and condition (\ref{conds5}) that 
if $2\le s\le r$,  then 
for $\theta\in\Theta^{\l,s}$, one has
$\theta_s=0$, i.e., $\theta$
is determined by $\theta_{\II{1}{s-1}}$ and $\theta_{\II{s+1}{r-1}}$, namely,
$\theta$ has $C_sC_{r-s}$ choices by the inductive assumption.
 Thus
$$
\#\Theta^\l=\mbox{$\sum\limits_{s=1}^{r+1}$}\#\Theta^{\l,s}
=C_{r-1}+\mbox{$\sum\limits_{s=2}^r$}C_sC_{r-s}+C_r=C_{r+1},
$$
where the last equality is a known
combinatorial identity, whose proof is omitted.
\end{proof}
The main result of this paper is the following theorem.
\begin{theorem}\label{theo-main}
Let $\l\in\Z_+^{m|n}$ be an $r$-fold atypical dominant integral weight. 
Then 
\equa{theo-main-1}
{
\mbox{the set of primitive weights of }\kac{\l}\ =\
\{L'_\theta(\l)\,\,|\,\,\theta\in\Theta^\l\}.
}
\end{theorem}
We shall prove Theorem \ref{theo-main} in Section \ref{proof}.
Let us look at the following
example.
\begin{example}\label{exam3.3}
Suppose again $\l$ is as in $(\ref{exam1})$. It is hardly possible to use
Theorem \ref{theo3.1} to determine the
primitive weights of $\kac{\l}$. However this job can be easily done by Theorem 
\ref{theo-main}.
Using Example
\ref{exam3.2}, $(\ref{exam-k-i})$ and $(\ref{R'-theta-f}),$ the primitive weights
of $\kac{\l}$ are the following $14$ weights:
the first $7$ weights are $($cf.~$(\ref{l-mu}))$
\equa{exam2}
{
\begin{array}{cc}
\!\!\!\!\!\!
(\ul{15},11,\ul{10},\ul{7},6,4,\ul{3}\,\,||\,\,\ul{3},5,\ul{7},8,\ul{10},\ul{15}),\!\!\!&
(\ul{15},11,\ul{10},\ul{7},6,4,\ul{2}\,\,||\,\,\ul{2},5,\ul{7},8,\ul{10},\ul{15}),\\[6pt]
\!\!\!\!\!\!
(\ul{15},11,\ul{10},6,4,\ul{2},\ul{1}\,\,||\,\,\ul{1},\ul{2},5,8,\ul{10},\ul{15}),\!\!\!&
(\ul{15},11,\,\ul{9},\,\ul{7},6,4,\ul{3}\,\,||\,\,\ul{3},5,\ul{7},8,\,\ul{9},\,\ul{15}),\\[6pt]
\!\!\!\!\!\!
(\ul{15},11,\,\ul{9},\,\ul{7},6,4,\ul{2}\,\,||\,\,\ul{2},5,\ul{7},8,\,\ul{9},\,\ul{15}),\!\!\!&
(\ul{15},11,\,\ul{9},\,6,4,\ul{2},\ul{1}\,\,||\,\,\ul{1},\ul{2},5,8,\,\ul{9},\,\ul{15}),\\[6pt]
\!\!\!\!\!\!
(\ul{15},11,\,\ul{7},\,6,4,\ul{2},\ul{1}\,\,||\,\,\ul{1},\ul{2},5,\ul{7},\,8,\,\ul{15}),\!\!\!
\end{array}\!\!\!\!\!
}
and the other $7$ weights are obtained from $(\ref{exam2})$
by changing $15$ to $14$ in all positions.
\end{example}
\begin{remark}\label{rema4.5}
One may regard Theorem \ref{theo-main} as the ``converse'' of Theorem
\ref{theo3.1} in the sense that Theorem
\ref{theo3.1} computes $a_{\l,\mu}$ for a given $\mu$ while Theorem
\ref{theo-main} computes $a_{\l,\mu}$ for a given $\l$.
\end{remark}
\section{A conjecture of Hughes et al}
\label{conjecture}
The purpose of this section is prove a  conjecture 
of Hughes, King and van der Jeugt \cite{HKV}.
We shall briefly recall some notions, which will be used throughout the section.
For more details, we refer to \cite{HKV, VHKT} (see also \cite{Cr2, VHKT0, VZ, SHK, S}).
\subsection{Composite Young diagram}
Let $\l=(\l^1\,\,|\,\,\l^2)$ be an $r$-fold atypical
dominant integral weight, written {\bf in terms of
notation (\ref{weight0})}, where $\l^1=(\l'_1,\l'_2,...,\l'_m)$ and
$\l^2=(\l'_{\cp1},\l'_{\cp2},...,\l'_{\cp n})$.
By the statements after (\ref{d-weight}), we can assume
\equa{standard}
{
\l'_1\ge\l'_2\ge...\ge\l'_m \ge 0\ge 
\l'_{\cp1}\ge\l'_{\cp2}\ge...\ge\l'_{\cp n},
}
Furthermore we can suppose $\l'_m$ and $-\l'_{\cp1}$ are large enough in order to
be able to perform boundary strip removals (see Subsection \ref{removal}).
Since both $\l^1$ and $-(\l^2)^R:=
(-\l_{\cp n}$, $-\l_{\cp n-1},...,\l_{\cp 1})$ are partitions,
this allows us to associate a {\it composite Young diagram} $F^{(\l^1|\l^2)}$
which is formed by joining the Young diagram $F^{\l^1}$ of $\l^1$ to the pointwise reflection
of the Young diagram $F^{-(\l^2)^R}$ of $-(\l^2)^R$. The part $F^{\l^1}$ is the
{\it covariant part} of $F^{(\l^1|\l^2)}$, and $F^{-(\l^2)^R}$ the {\it contravariant part}
of $F^{(\l^1|\l^2)}$. The diagram $F^{(\l^1|\l^2)}$ is 
{\it standard} if (\ref{standard}) holds.
\begin{example}\label{exam4.1-0}
Suppose $\l$ is as in $(\ref{exam1}),$ or $\l=(8,5,5,3,3,2,2{\ssc\,}|{\ssc\,}-2,
{\ssc\!}-3,{\ssc\!}-4,{\ssc\!}-4,{\ssc\!}-5,{\ssc\!}-9)$
in terms of notation $(\ref{weight0}),$ then the composite Young diagram of $\l$ is
\\[4pt]\cl{
$
\begin{array}{rl}
-(\l^2)^R \ \ \ =\ \ \
\begin{array}{r}
\Box
\vs{-6pt}\\
\Box
\vs{-6pt}\\
\Box
\vs{-6pt}\\
\Box
\vs{-6pt}\\
\put(3,1.8){$\sc3$}\Box\!\Box
\vs{-6pt}\\
\put(3,1.8){$\sc3$}\Box\!\put(1,1.8){$\sc3$}\!\Box\!\put(1,1.8){$\sc3$}\!\Box\!\Box
\vs{-6pt}\\
\put(3,1.8){$\sc3$}\Box\!\put(1,1.8){$\sc3$}\!\Box\!\Box\!\Box\!\Box
\vs{-6pt}\\
\put(3,1.8){$\sc1$}\Box\!\put(1,1.8){$\sc3$}\!\Box\!\Box\!\Box\!\Box\!\Box
\vs{-6pt}\\
\put(3,1.8){$\sc3$}\Box\!\put(1,1.8){$\sc3$}\!\Box\!\Box\!\Box\!\Box\!\Box
\end{array}
\\\!\!\!\!&\!\!\!\!\!\!\!\!
\begin{array}{l}
\Box\!\Box\!\Box\!\Box\!\Box\!\Box\!\Box\!\Box
\vs{-6pt}\\
\Box\!\Box\!\Box\!\Box\!\Box
\vs{-6pt}\\
\Box\!\Box\!\Box\!\put(1,1.8){$\sc3$}\!\Box\!\put(1,1.8){$\sc3$}\!\Box
\vs{-6pt}\\
\Box\!\Box\!\put(1,1.8){$\sc3$}\!\Box
\vs{-6pt}\\
\Box\!\put(1,1.8){$\sc3$}\!\Box\!\put(1,1.8){$\sc3$}\!\Box
\vs{-6pt}\\
\put(3,1.8){$\sc3$}\Box\!\put(1,1.8){$\sc3$}\!\Box
\vs{-6pt}\\
\put(3,1.8){$\sc3$}\Box\!\put(1,1.8){$\sc1$}\!\Box
\end{array}
\ \ \ =\ \ \ \l^1\ ,
\end{array}
$}\\[4pt]
where the labeled boxes will be explained in Subsection \ref{removal}.
\end{example}
\subsection{Permissible code}
To determine composition factors of Kac modules, Hughes et al \cite{HKV}
introduced the notion
of permissible codes, which we recall below.
\begin{definition}\label{code}
Suppose $\l$ is an $r$-fold atypical dominant integral weight.
A permissible code $\mu_c$ for $\l$ is an array of length $r$, each element of the
array consisting of a non-empty column of increasing labels taken from
$\{0,1,...,r\}$. The first element of a
column is called the top label. A permissible
code $\mu_c$ must satisfy the rules:

{\rm(i)}
The top label of column $s$ can be $0,s$ or $a$ with
$s<a$; the first case can occur only if column $s$ is zero, while
the last case can occur only if $c_{s,t}=\qq$ with $a$
the top label of column $t$ for some $t>s$.

{\rm(ii)}
Let $s<t,\,c_{s,t}=c_{s+1,t}=...=c_{t-1,t}=\cc$. If the top label of column $t$
is $a$ with $t\le a$, then $a$ must appear somewhere below the top entry
of column $s$.

{\rm(iii)}
If $s$ appears in any column then the only labels
which can appear below $s$ in the same column are those $t$ with $s<t$,
for which $t$ is the top label of column $t$ and
$c_{s,t}= \cc$.

{\rm(iv)}
If the label $s$ appears in more than one column and $t$
appears immediately below $s$ in one such column, then it must do so in all
columns containing $s$.

{\rm(v)}
Let $s<t<u$ and
$c_{s,t} = \qq,\,c_{t,u} = \qq$
{\rm(}so, $c_{s,u}=\qq${\rm)}. If the top label of
column $s$ is the same as that of column $u$ and it is nonzero
then the top label of column $t$ is not $0$.

{\rm(vi)}
Let $s < t < u < v$ with top labels $a,b,a,b$
respectively, $a\ne0\ne b$. If $a < b$ then columns $s$ and $u$ must
contain $b$; if $a > b$ then columns $t$ and $v$ must contain $a$.

{\rm(vii)}
If a column has two nonzero labels, then the last label of this column must
appear in the next column.
\end{definition}
\begin{example}\label{exam4.1}
Suppose $\l$ is as in $(\ref{exam1})$.
Using the rules in Definition \ref{code},
we find the following $14$ permissible codes $\mu_c$:
\equa{Ex1-code}{
\begin{array}{lllllll}
^{\dis0\ 0\ 0\ 0}\ &\ ^{\dis1\ 0\ 0\ 0}\ &\ ^{\dis1\ 2\ 0\ 0}_{_{\dis2}}\ &\ 
^{\dis0\ 0\ 3\ 0}\ &\ ^{\dis1\ 0\ 3\ 0}\ &\ ^{\dis1\ 2\ 3\ 0}_{_{\dis2}}\ &\
^{\dis1\ 3\ 3\ 0}_{_{\dis3}}
\\[12pt]
^{\dis0\ 0\ 0\ 4}\ &\ ^{\dis1\ 0\ 0\ 4}\ &\ ^{\dis1\ 2\ 0\ 4}_{_{\dis2}}\ &\ 
^{\dis0\ 0\ 3\ 4}\ &\ ^{\dis1\ 0\ 3\ 4}\ &\ ^{\dis1\ 2\ 3\ 4}_{_{\dis2}}\ &\
^{\dis1\ 3\ 3\ 4}_{_{\dis3}}
\end{array}
}
\end{example}
\begin{remark}\label{rmk-code}
{\rm(1)} Rule {\rm(vii)} in Definition {\rm\ref{code}} does not appear 
in the definition of a code in \cite{HKV}. 
However we observe the fact that in all examples of determining codes given
in \cite{HKV},
this rule was implicitly applied.
We realize that if this rule is not included in the definition, one would
produce two ``codes'' such as 
$({\ssc\,}^{1\ 0\ 3\ 0}_{3}{\ssc\,})$ and
$({\ssc\,}^{1\ 3\ 3\ 0}_{3}{\ssc\,})$ $($cf.~$(\ref{Ex1-code}))$
which would correspond to the same sequence of boundary strip removals $($see Subsection
\ref{removal}$)$.

{\rm(2)} From the proof of Theorem \ref{theo-main}, we shall see that a code 
must satisfy the
following rule stronger than rule {\rm(vii):}
\begin{itemize}
\item[{\rm(vii)$'$}]
If a column has two nonzero labels, then all labels, except possibly the top label,
of this column must appear in the next column.
\end{itemize}
\end{remark}
\begin{remark}\label{rmk-code-0}
Rules {\rm(iv)}, {\rm(vi)} and {\rm(vii)$'$} imply the following rule
stronger than rule {\rm(vi):}
\begin{itemize}
\item[{\rm(vi)$'$}]
Let $s < t < u < v$ such that columns $s,\,t,\,u,\,v$ contain labels $a,\,b,\,a,\,b$
respectively with $a\ne0\ne b$. If $a < b$ then columns $s$ and $u$ must
contain $b$; if $a > b$ then columns $t$ and $v$ must contain $a$.
\end{itemize}
\end{remark}
\subsection{Boundary strip removal}
\label{removal}
For a code $\mu_c$ for $\l$, it corresponds to a weight $\mu$ defined as follows.

First $\mu_c$ corresponds to a sequence of 
{\it coordinated boundary strip removals} \cite{HKV} defined by:
For each $s=1,2,...,r$, if label $s$ appears in $\mu_c$, 
first we set
\equa{a-s-l}
{a_s=\min\{p\in\II{1}{r}\,\,|\,\,\,\mbox{label }s\mbox{ appears in column }p\},
}
(i.e., column $a_s$ is the first column which contains label $s$),
and 
we use notation (\ref{aty-roots}),
then 
a {\it coordinated boundary strip removal starting from the $s$-th atypical $\g_s$ and ending at 
the $a_s$-th atypical $\g_{a_s}$} is performed, that is,
two
boundary strip removals are simultaneously performed
on $F^{\l^1}$ and on $F^{-(\l^2)^R}$
which start respectively from the $m_s$-th row of $F^{\l^1}$ and 
the $(\cp n+1-n_s)$-th column of $F^{-(\l^2)^R}$ and continue
until they pass the $m_{a_s}$-th row  of $F^{\l^1}$ and 
the $(\cp n+1-n_{a_s})$-th column of $F^{-(\l^2)^R}$, then
continue until the remaining composite Young diagram is standard.
Then $\mu$ is defined to be the weight whose composite Young diagram is the remaining diagram.

Thus we obtain a correspondence
\equa{corres}
{\mu_c\mapsto\mu.
}
\begin{example}\label{exam4.2}
If $\l$ is as in Example \ref{exam4.1-0}, and $\mu_c$ is the
$7$-th code $({\ssc\,}^{1\ 3\ 3\ 0}_3{\ssc\,})$
in $(\ref{Ex1-code})$, then the boxes in the
sequence of coordinated boundary strip removals are labeled in the
diagram $F^{(\l^1|\l^2)}$ of  Example \ref{exam4.1-0}. Thus
\equa{exam-mu}
{\mu=(8,5,3,2,1,0,0\,\,|\,\,0,0,-2,-3,-3,-9).
}
which is the same as the last weight in $(\ref{exam2}).$
\end{example}
Similarly, an element $\theta\in\Theta^\l$ also corresponds to
a sequence of coordinated boundary strip removals, such that for
each $s=1,2,...,r$, if $\theta_s\ne0$, then
a coordinated boundary strip removal starting from $\g_s$ and ending at 
$\g_{s+1-\theta_s}$ is performed.
Then we see that the remaining diagram
is the composite Young diagram of $\mu$, where $\mu=L'_\theta(\l)$.
\begin{example}\label{exam4.3}
If $\theta=(1,0,3,0)$ as in $(\ref{exam2-1})$, then again we have the
sequence of coordinated boundary strip removals in Example \ref{exam4.1-0},
and $\mu$ is as in $(\ref{exam-mu})$.
\end{example}
\begin{remark}\label{addition}
Parallel to boundary strip removals, one can introduce 
the notion of boundary strip additions $($we omit the precise
definition here
since we shall not use this later$)$, 
such that for all $r$-fold atypical
$\mu\in\Z_+^{m|n}$, each $\theta'\in\{0,1\}^r$ ``corresponds'' 
$($under some regulations$)$ to 
a sequence of coordinated boundary strip additions performed on the 
composite Young diagram of $\mu$. The resulting diagram is the diagram of
$\l$, where $\l=R_{\theta'}(\mu)$ $($cf.~Theorem \ref{theo3.1}$)$.
For example, if $\mu$ is the weight corresponding to the un-labeled boxes
in Example \ref{exam4.1-0} $($i.e., $\mu$ is the weight in $(\ref{exam-mu}))$, 
and
$\theta'=(1,1,0,0)$, then $\l=R_{\theta'}(\mu)$ is the weight in the example,
where labels $3$ and $1$ in the diagram shall be changed to $1$ and $2$ respectively
in order that the labeled boxes correspond to the boundary strip additions
$($cf.~Lemma \ref{ll-1} below$)$.
\end{remark}
\subsection{A conjecture of Hughes et al}
As an application of Theorem \ref{theo-main}, we prove the
following theorem which was a conjecture
put forward by Hughes, King and van der Jeugt
in \cite{HKV} as the result of in depth research carried out by the authors over several years time.
\begin{theorem}\label{code-primitive}
The correspondence
$(\ref{corres})$ is a $1$ - $1$ correspondence 
between the set of primitive weights of $\kac{\l}$ and
the set of permissible codes for $\l$.
\end{theorem}
Partial result of this theorem was obtained in \cite{SHK} (also cf.~\cite{S}),
where it was proved that an ``unlinked code'' (i.e., a code that does not
have two columns with the same nonzero top label) 
corresponds to a ``strongly'' primitive weight
(i.e., a primitive weight
whose primitive vector is a highest weight vector in the Kac-module); 
moreover the primitive vector
corresponding to the unlinked code is precisely constructed.

\vskip6pt\ni{\it Proof of Theorem \ref{code-primitive}}.
Let
\equa{code-set}
{
\cC= \mbox{ the set of permissible codes for }\l.
}
We  establish a $1$ - $1$ correspondence between $\cC$ and $\Theta^\l$
(then Theorem \ref{code-primitive} follows from  Theorem \ref{theo-main}).
Let $\mu_c\in\cC$ be a code. We define $\theta\in\Theta^\l$ as follows:
for $s=1,2,...,r$, 
\begin{eqnarray}
\label{theta-not-0}
\theta_s\ne0&\Lra&
\mbox{label $s$ appears in the code $\mu_c$, and
in this case}
\\
\label{theta-s-code}
&&
\theta_s=s+1-a_s,
\end{eqnarray}
where $a_s$ defined in (\ref{a-s-l}) is the first number whose column contains label
$s$ (and obviously column $s$ is the last number whose column
contains label $s$ as the top label:
note that if a label $s$ appears in a code $\mu_c$, it must be the top label of column $s$). 
For example, if $\mu_c$ is a code in (\ref{Ex1-code}), then $\theta$ is a
corresponding element in (\ref{exam2-1}).

We want to prove that $\theta$ is indeed in $\Theta^\l$ by verifying each condition
of (\ref{conds1})--(\ref{conds5}). This will be done by several claims.
So suppose $\theta_s\ne0$ for some $s\in\II{1}{r}$.
\VB
{\bf Claim 1}. The second condition of (\ref{conds1}) holds, namely, $c_{a_s,s}\ne\nn$.
\VE
Since label $s$ appears in column $a_s$,
if it is the top label, by rule
(i) and the fact that $\qq$-relationship is transitive we have $c_{a_s,s}=\qq$.
Otherwise, let $p$ be the top label of column $a_s$, then $p<s$ and $c_{p,s}=\cc$
by rule (iii). But rule (i) says that $c_{a_s,p}=\qq$, which together with the relation
$c_{p,s}=\cc$ implies that $c_{a_s,s}=\cc$. This proves Claim 1.
\VB
{\bf Claim 2}. The first condition of (\ref{conds1}) holds, namely, $c_{a_s-1,s}\ne\cc$.
\VE
Suppose conversely $c_{a_s-1,s}=\cc$. 
\VB
{\it Case $($a$)$: First suppose $c_{a_s-1,a_s}=\cc$}.
Since label $s$ appears in column $a_s$ of $\mu_c$,
by rules (ii) and (iv), all labels of column $a_s$ must appear in column $a_s-1$,
in particular $s$ appears in column $a_s-1$. This 
in turn contradicts definition (\ref{a-s-l}). Thus this case does not occur.
\VB
{\it Case $($b$)$: Suppose  $c_{a_s-1,a_s}\ne \cc$}.
If $c_{a_s,s}=\qq$, then this and the relation $c_{a_s-1,a_s}\ne \cc$ imply that
$a_{a_s-1,s}\ne \cc$ by definition (\ref{c'-st}), a contradiction with the assumption. Thus 
$c_{a_s,s}=\cc$ by Claim 1.

If $c_{p,s}=\cc$ for all $p\in\II{a_s}{s-1}$, then by rule (ii), label $s$ must appear in
column $a_p-1$, a contradiction with the definition (\ref{a-s-l}).
So $c_{p',s}\ne \cc$ for some $p'\in\II{a_s}{s-1}$.
Set 
\equa{ppp}
{
p'=\min\{p'\in\II{a_s}{s-1}\,\,|\,\,c_{p',s}\ne \cc\}
}
to be the smallest whose corresponding atypical root $\g_{p'}$  of $\l$
is not $\cc$-related to 
the $s$-th atypical root $\g_s$.
Then $c_{p',s}=\qq$ (otherwise $c_{p'-1,s}$ cannot be $\cc$, contradicting (\ref{ppp})).
From this and definition (\ref{ppp}), we have
\equa{CCCCCCCCC}
{c_{p,p'}=c_{p,s}=\cc\mbox{ \ for all \ }p\in\II{a_s-1}{p'-1}.
} 
So rules (ii) and (iv) show that 
\equa{SHOW}
{\mbox{all nonzero labels
in column $p'$ must appear in columns $a_s$ and $a_s-1$.
}}
Thus the following subclaim means
that $s$ appears in column $a_s-1$, this contradiction with 
definition (\ref{a-s-l}) implies Claim 2.
\VB
{\bf Subclaim 2a)}.
The top label of column $p'$ is $s$.
\VE
To prove this subclaim, suppose $u$ is the last label of column $a_s$.
First assume that $u\ne s$. Then $s<u$ and $c_{s,u}=\cc$ by rule (iii). This, together with
the relation $c_{p',s}=\qq$ and (\ref{CCCCCCCCC}),
shows that
$c_{p,u}=\cc$ for $p\in\II{a_s-1}{p'}$, which in turn implies that
$u$ appears in column $p$ by rule (vii). Thus $u$ appears in column $p'$ and it is
not the top label of $p'$ by rule (i).
So let $u'$ be the top label of $p'$. Then $c_{p',u'}=\qq$ and so $c_{s,u'}=\qq$ (if $s\le u'$)
or $c_{u',s}=\qq$ (if $u'<s$) by the fact that $c_{p',s}=\qq$.
By rule (iii), two labels whose corresponding atypical roots are $\qq$-related cannot
appear in the same column, but both $s$ and $u'$ appear in column $a_s$ by (\ref{SHOW}).
We have $u'=s$, and thus obtain the subclaim in this case.

Next assume $s=u$ is the last label of coumn $a_s$. Then rule (vii) implies that
$s$ appears in column $p$ for all $p\in\II{a_s-1}{p'}$. In particular, $s$
must be the top label of column $p'$. The subclaim is proved.

(The proof of the subclaim also shows that a code satisfies rule
(vii)$'$ in Remark \ref{rmk-code}.)
\VB
{\bf Claim 3}. Condition (\ref{conds2}) holds.
\VE
First suppose $\theta_p>\theta_s-s+p$
 for some $p\in\II{a_s}{s{\ssc\!}-{\ssc\!}1}$.
Then we have $a_p{\ssc\!}<{\ssc\!}a_s{\ssc\!}<{\ssc\!}p{\ssc\!}<{\ssc\!}s$
 and columns $a_p,\,a_s,\,p,\,s$ 
contain labels $p,\,s,\,p,\,s$ respectively. 
By rule (vi)$'$ in Remark \ref{rmk-code-0}, $s$ appears in column $a_p$ (which is $<a_s$), a contradiction with 
definition (\ref{a-s-l}). Next suppose
$\theta_p=\theta_s-s+p$. Then $a_p=a_s$ and so both $p$ and $s$ appear in column $a_s$.
By rule (iii), we must have $c_{p,s}=\cc$. This proves Claim 3.
\VB
{\bf Claim 4}. Condition (\ref{conds4}) holds.
\VE
Assume conversely for some $p\in\II{a_s}{s-1}$,
$c_{p,s}\ne \nn$, but
\equa{theta-p-0}
{\begin{array}{l}
\theta_p=0\mbox{ and no }p'\in\II{p+1}{s}
\mbox{ with }c_{p,p'}=\qq\mbox{ and }\theta_{p'}\ge p'+1-p,
\\[4pt]
\mbox{and }p\ne a_s\mbox{ or }c_{a_s,s}\ne \qq
\mbox{ \ \ (cf.~condition (\ref{conds5}))}.
\end{array}
}
Condition (\ref{theta-p-0}) implies 
\equa{zero}
{
\mbox{column $p$ of $\mu_c$}\ =\ 0.
}
Claim 1 allows us to consider the following two cases: 
\VB
{\it Case $($a$)$: First suppose $c_{a_s,s}=\qq$.}
Then rule (v) and (\ref{zero}) mean that $c_{p,s}\ne \qq$, i.e., $c_{p,s}=\cc$.
This together with rule (ii) and (\ref{zero}) implies
that
(cf.~the arguments in Case (b) of Claim 2)
 there exists $p_1\in\II{p+1}{s-1}$ such that $c_{p_1,s}\ne \cc$.
Let
\equa{p-1-}
{
p_1=\min\{p_1\in\II{p}{s-1}\,\,|\,\,c_{p_1,s}\ne \cc\}
}
be the smallest whose corresponding atypical root $\g_{p_1}$ of $\l$
is not $\cc$-related to $\g_s$.
If $c_{p_1,s}=\nn$, then from definition (\ref{c'-st}),
it follows that $c_{p_1-1,s}\ne \cc$,
contradicting (\ref{p-1-}).
Thus $c_{p_1,s}=\qq$.
Using this and definition (\ref{p-1-}), we have 
$c_{p',p_1}=c_{p',s}=\cc$ for all $p'\in\II{p}{p_1-1}$. This, 
together with the fact that column $p_1$ is nonzero 
(which is derived by rule (v) and the relation $c_{p_1,s}=\qq$),
contradicts (\ref{zero}) by rule (ii).
\VB
{\it Case $($b$)$: Suppose $c_{a_s,s}=\cc$.}
This time instead of defining $p_1$ in (\ref{p-1-}), we define $p'_1$ by
\equa{p-1-1}
{
p'_1=\min\{p'_1\in\II{a_s}{s-1}\,\,|\,\,c_{p'_1,s}\ne \cc\}.
}
The fact (\ref{zero}) ensures that such $p'_1$ exists 
by using rule (ii) since column $s$ is nonzero,
and $c_{p'_1,s}=\qq$ as in Case (a).
Now rules (ii), (vii) and rule (vi)$'$ in Remark \ref{rmk-code-0} show that
$s$ appears in column $p'$ for all $p'\in\II{a_s}{p'_1}$. In particular $p'_1<p$ by 
(\ref{zero}).
Now by the same arguments in Case (a) (with $a_s$ replaced with $p'_1$), one obtains
a contradiction. This proves Claim 4.
\VB
{\bf Claim 5}. Condition (\ref{conds5}) holds.
\VE
To prove this, suppose $a_s<s$ and $c_{a_s,s}=\qq$.
Since $s$ appears in column $a_p$
(cf.~the first statement after Claim 1), by rules (i) and
(iii), $s$ must be the top label of column $a_s$. 
This in particular implies that label $a_s$ cannot appear anywhere
in $\mu_c$ (cf.~the first statement after (\ref{theta-s-code})), 
i.e., $\theta_p=0$ for $p=a_s$. Claim 5 is proved.
\VB
Denote the obtained $\theta$ by $\theta_{\mu_c}$. The above claims show that we 
have a map
\equa{map-code}
{
\cC\to\Theta^\l:\ \mu_c\mapsto\theta_{\mu_c}.
}

Conversely, suppose $\theta\in\Theta^\l$. We define a code $\mu_c\in\cC$ as follows.

For $s=1,2,...,r$, if $\theta_s\ne0$, 
this time we first set  $a_s$ to be $s+1-\theta_s$ (cf.~(\ref{theta-s-code})),
then the code $\mu_c$ should satisfy:
label $s$ must appear in column $s$ (as the top
label) of $\mu_c$, and appear in column $a_s$ (columns $a_s$
and $s$ are respectively the first and last columns where label $s$ appears);
and 
for all $p\in\II{a_s+1}{s-1}$, label $s$ appears in column $p$ if and only if
$c_{p,s}=\cc$, or $c_{p,s}=\qq$ but there does not exist a smaller label 
$s'\in\II{p}{s-1}$ which appears in this column.

It is a little tedious but straightforward routine to 
check that the above uniquely defines a code $\mu_c$, denoted by 
$\mu_c^\theta$,
satisfying rules (i)--(vii), and that the correspondence 
\equa{another-cor}
{\theta\mapsto\mu_c^\theta}
obtained
in this way is the inverse of the map (\ref{map-code}).
We omit the details since they are mainly the reverse of the above arguments.
\hfill$\Box$ 
\section{Proof of the main theorem}
\label{proof}
\def\wh#1{\ul{#1}{\sc\,}}
The aim of this section is to give a proof of Theorem \ref{theo-main}.
This will be done by several lemmas.
Having given a proof of Theorem \ref{code-primitive}
helps us in understanding the arguments below.

Parallel to definition (\ref{k--i}),
we introduce $r$-tuples $(\wh k_1^{(\nu)},\wh k_2^{(\nu)},...,
\wh k_r^{(\nu)}),\,\nu\ge1$, of positive integers defined by
\equa{k---i}
{
\wh k^{(\nu)}_s=\min\left\{k>0\,\left|\ \#\left({}^{^{^{^{\,\!}}}}\II{\l_{m_s}}{\l_{m_s}+k}\bs S(\l)\right)
=\nu\right\}.\right.
}
Thus $\l_{m_s}+\wh k^{(\nu)}_s$ is the $\nu$-th smallest integer bigger than $\l_{m_s}$ and not in $S(\l)$,
and the definition of $k_s$ in (\ref{k-i}) implies
\equa{k-i-i}
{
k_s=\wh k_s^{(p_s+1-s)}\mbox{ for }s=1,2,...,r.
}

By Theorem \ref{theo3.1}, the proof of (\ref{theo-main-1}) is equivalent to proving
that
for  $r$-fold atypical weights $\l,\mu\in\Z_+^{m|n}$,
\equa{to-prove}
{
\mu=L'_\theta(\l)\mbox{ for some }\theta\in\Theta^\l
\ \Lra\
\l=R_{\theta'}(\mu)\mbox{ for some }\theta'\in\{0,1\}^r.
}
\begin{convention}\label{con5}
For convenience,
unless it is specified, we shall always use 
the same notation with a tilde to denote any element or concept
associated with $\mu$; for instance, 
$\wt k_i=k_i^\mu$ $($recall $(\ref{k-i})$ and Convention \ref{conv}$)$, $\wt p_s=p^\mu_s$,
$\wt m_s=m^\mu_s$, etc.
\end{convention}
First suppose
$\mu=L'_\theta(\l)$ for some $\theta\in\Theta^\l$.
For $s=1,2,...,r$, we define
\begin{eqnarray}
\label{def-pi-i}
&&
a_s=s+1-\theta_s,
\\
\label{def-N-s}
&&
N_{s,p}=\#\{p'\in\II{a_s}{p-1}\,\,|\,\,\theta_{p'}\ne0,\,a_{p'}=a_s\}\mbox{ if }\theta_s\ne0,
\end{eqnarray}
for $p\in\II{a_s}{s-1}$,
and set $N_{s,p}=-1$ if $\theta_s=0$. We simply denote $N_s=N_{s,s}$.
In terms of the code $\mu_c^\theta$ defined in
(\ref{another-cor}),
$a_s$ defined here is
the first number whose column contains label $s$ in the code $\mu_c^\theta$  (cf.~(\ref{a-s-l})),
and
$N_{s,p}$ is the number of those labels $p'$ which is smaller than $p$ and
which first appear in column $a_s$ (``first appearance'' means ``not appear in a smaller
column'').
\begin{lemma}\label{l5.1}
If $\theta_s\ne0$, then $($recall definition $\ell_{s,t}$ in $(\ref{c'-st}))$
\begin{eqnarray}
\label{l5.1-1}
&&
a_s\le a_p\le p\mbox{ \ for all \ }p\in\II{a_s}{s-1}\mbox{ with }\theta_p\ne0,
\\
\label{l5.1-2}
&&
N_s=s-a_s-\ell_{a_s,s}\mbox{ \ and}
\\
\label{l5.1-3}
&&
N_{s,p}<N_s\mbox{ if }p\in\II{a_s}{s-1}.
\end{eqnarray}
\end{lemma}
\begin{proof}
The first equation follows from condition (\ref{conds2}) and the last from definition
(\ref{def-N-s}).
The second equation can be proved by induction on $s$ as follows.

Set $N'_s=s-a_s-\ell_{a_s,s}$.
If $N'_s=0$, i.e., $c_{a_s,s}=\qq$, then (\ref{l5.1-2}) follows from condition
(\ref{conds5}). Suppose $N'_s>0$, i.e., $c_{a_s,s}=\cc$. Condition (\ref{conds4}) shows that
there exists $p'\in\II{a_s}{s-1}$ with $\theta_{p'}\ne0$ such that $a_{p'}=a_s$. Let $p'$
be the maximal such number. This definition of $p'$ implies that
$N_{s,p'}=N_s-1$. But in fact we have $N_{s,p'}=N_{p'}$ from the definition 
(\ref{def-N-s}).
Thus $N_{p'}=N_s-1$. Since $p'<s$, by inductive assumption on $p'$, 
we can suppose $N'_{p'}=N_{p'}$, i.e,
$p'-a_s-\ell_{a_s,p'}=N_s-1$. Thus the proof of (\ref{l5.1-2}) is reduced to 
proving
\equa{REMAIN}
{
p'-\ell_{a_s,p'}=
s-\ell_{a_s,s}-1.
}
By condition (\ref{conds2}), 
 \equa{CAC}
{c_{p',s}=\cc\mbox{, i.e., }\ell_{p',s}<s-p'.
}
By the maximal choice of $p'$,
we must have $\ell_{p',s}=s-p'-1$
(the arguments to prove this are similar to those given in Case (b) of Claim 2 in
the proof of Theorem \ref{code-primitive}, thus omitted), 
which is equivalent to (\ref{REMAIN}).
\end{proof}

We define
\begin{eqnarray}
\label{pi-i}&&
\pi_s=\#\{p\in\II{s}{r}\,\,|\,\,s=a_p=p+1-\theta_p\},
\\
\label{define-t'}
&&
\theta'_s=\biggl\{\begin{array}{ll}
1&\mbox{if }\exists\,p\in\II{1}{s}\mbox{ such that }
p\le s< p+\pi_p,
\\[2pt]
0&\mbox{otherwise}.
\end{array}
\end{eqnarray}
Here $\pi_s$ is the number of those labels $p$ which first appear in column $s$ of the code
$\mu_c^\theta$ 
(thus in particular $N_s=\pi_{a_s}$).
Definition (\ref{define-t'}) means that
$\theta'_s\ne0$ if and only if 
column $s$ is bigger than but ``close'' to 
a (unique) column $p$ (where some labels first appear), here
``close'' means that the distinct $s-p$ is smaller than the number $\pi_p$ of 
those labels first appearing in column $p$.

For example, if $\theta=(1,0,3,0)$ in (\ref{exam2-1}), then
$\mu^\theta_c=({\ssc\,}^{1\ 3\ 3\ 0}_3{\ssc\,})$ and
\equa{theta'-exam}
{
(\pi_1,\pi_2,\pi_3,\pi_4)=(2,0,0,0), \mbox{ \ and so \ }
\theta'=(1,1,0,0).}
Note from definitions (\ref{def-pi-i}) and (\ref{define-t'}) that
\equa{equal-num}
{
\mbox{$\sum\limits_{s\in\II{1}{r}}$}\pi_s=
\#\{s\in\II{1}{r}\,\,|\,\,\theta_s\ne0\}=
\#\{s\in\II{1}{r}\,\,|\,\,\theta'_s\ne0\}=:N,
}
which is the number of nonzero $\theta_s$'s (or nonzero $\theta'_s$'s).
Lemma \ref{l5.1} allows us to 
divide the index set $\II{1}{r}$ into
some subsets: $\II{1}{r}=\VV_1\cup\VV_2$ and
\equa{subset}
{
\VV_2=\II{\ell'_1}{\ell_1}\cup\II{\ell'_2}{\ell_2}
\cup\cdots\cup\II{\ell'_\nu}{\ell_\nu}
}
for some $\nu$ and some $\ell'_1\le \ell_1<\ell'_2\le\ell_2<...<
\ell'_\nu\le \ell_\nu$,
such that $\theta_p=0$ for $p\in\VV_1$ and
\equa{either-or}
{
\theta_{\ell_i}\ne0\mbox{ and }
\ell_i=
\max\{p\in\II{\ell'_i}{r}\,\,|\,\,a_p=\ell'_i\},\,i=1,...,\nu.
}
(The above means that $\ell_i$ is the largest label which first
appears in column $\ell'_i$ of the code $\mu_c^\theta$, and 
this in turn means that all labels of column
$\ell'_i$ first appear in this column.)
\begin{lemma}\label{lemm5.1}
$\l=R_{\theta'}(\mu).$
\end{lemma}
\begin{proof}
\def\Eta{\varphi}
Denote $\Eta=R_{\theta'}(\mu)$.
By (\ref{l-mu}), it suffices to prove
$T(\l)=T(\Eta)$
(cf.~definition (\ref{set-R(f)})), 
which in turn is equivalent to
\equa{suffice-prove}
{T(\Eta)\subset T(\l)
}
by the dominance of $\Eta$.
Equations (\ref{Main-t}) and (\ref{Main--t}) imply  (cf.~Convention \ref{con5})
\equa{TTT}
{
T(\Eta)=\{\mu_{\wt m_s}+\theta'_s\wt k_s\,\,|\,\,s\in\II{1}{r}\}
\mbox{ \ and \ }
T(\mu)=\{\l_{m_s}-\low k\!\!\!k\!\!\!k^{(\theta_s)}_s\,\,|\,\,s\in\II{1}{r}\},
}
where $\low k\!\!\!k\!\!\!k^{(\theta_s)}_s$ is defined after (\ref{set-l-s-}).
We want to prove (\ref{suffice-prove}).
Nothing needs to be done if $N=0$. So suppose $N>0$.
We first need to list the order of integers in $T(\mu)$ in order to determine
$\mu_{\wt m_s}$ for $s\in\II{1}{r}$. This is done by the following claims.
\VB
{\bf Claim 1}. Let $p\in\II{\ell'_i}{\ell_i},\,s\in\II{\ell'_j}{\ell_j}$ for some
$i<j$. Then $\l_{m_p}<\l_{m_s}-\low k\!\!\!k\!\!\!k^{(\theta_s)}_s$.
\VE
If $\theta_s=0$ there is nothing to do since $\l_{m_p}<\l_{m_s}$. Suppose $\theta_s>0$.
Then $p\le\ell_i<\ell'_j\le a_s$ and $c_{a_s-1,s}\ne \cc$, i.e., $\ell_{a_s-1,s}\ge s-(a_s-1)
=\theta_s$. This and the definition of $\low k\!\!\!k\!\!\!k^{(\theta_s)}_s$ show
that $\l_{m_s}-\low k\!\!\!k\!\!\!k^{(\theta_s)}_s$
which is the $\theta_s$-th largest integer smaller than $\l_{m_s}=\l^{(s-1)}_{m_s}$ 
and not in $S(\l^{(s-1)})$ (recall the definition of $\l^{(s)}$ in (\ref{set-l-s-})),
is bigger than $\l_{m_p}$. This is the claim.
\VB
{\bf Claim 2}. 
Let $p\in\II{\ell'_i}{\ell_i-1}$ and $s=\ell_i$. 
Then $\l_{m_s}-\low k\!\!\!k\!\!\!k^{(\theta_s)}_s<\l_{m_p}-\low k\!\!\!k\!\!\!k^{(\theta_p)}_p$.
\VE
Suppose $a_p=a_s$. The proof of Lemma \ref{l5.1} shows that
$\ell_{p,s}<s-p$ (cf.~(\ref{CAC})). Thus $\theta_s-\theta_p=s-p>\ell_{p,s}$, 
which in turn
says that ``the $\theta_s$-th largest integer smaller than $\l^{(s-1)}_{m_s}$ and not in $S(\l^{(s-1)})$''
is smaller than ``the $\theta_p$-th largest integer smaller than $\l^{(p-1)}_{m_p}$ and not in $S(\l^{(p-1)})$''.
This is Claim 2 by the definition of $\low k\!\!\!k\!\!\!k^{(\nu)}_s$.

If $a_s<a_p$, as in the proof of Claim 1 one can show that $\l_{m_{a_s}}<\l_{m_p}-\low k\!\!\!k\!\!\!k^{(\theta_p)}_p$.
But obviously $\l_{m_s}-\low k\!\!\!k\!\!\!k^{(\theta_s)}_s<\l_{m_{a_s}}$. Claim 2 is proved.
\VB
{\bf Claim 3}.
For $p,s\in\II{\ell'_i}{\ell_i}$, 
$\l_{m_s}-\low k\!\!\!k\!\!\!k^{(\theta_s)}_s<\l_{m_p}-\low k\!\!\!k\!\!\!k^{(\theta_p)}_p$ if and only
$a_s<a_p$ or $a_s=a_p$ but $p<s$.
\VB
This claim follows from the same arguments in the proofs of the above two claims.
\VB
The above three claims completely determine the order of the elements in the set $T(\mu)$.
From this and the definition (\ref{define-t'}), one can see that
$\theta'_p=0$ if and only $\theta_s=0$ for some $s$ such that $\mu_{\wt m_p}=\l_{m_s}$.
Thus 
to prove (\ref{suffice-prove}), by (\ref{TTT}),
we only need to consider
elements $\mu_{\wt m_s}+\theta'_s\wt k_s$ of $T(\eta)$ with $\theta'_s\ne0$.
We denote 
\equa{T-i-mu}
{T_i(\mu)=\{\l_{m_s}-\low k\!\!\!k\!\!\!k^{(\theta_s)}_s\,\,|\,\,s\in\II{\ell'_i}{\ell_i}\}.
}
By using induction on $\nu$, we shall only need to consider the set 
$T_i(\mu)$ for a particular $i$. By restricting to 
$\II{\ell'_i}{\ell_i}$, we may
regard $\II{\ell'_i}{\ell_i}$ as the whole set $\II{1}{r}$,
i.e., $\ell'_i=1,\,\ell_i=r$, in order to simplify
notations. 
In such case, $\theta_r=r$,
and $\theta'_1=1$.
Then Claims 2 and 3 show that 
$\l_{m_r}-\low k\!\!\!k\!\!\!k^{(r)}_r$ 
(in this case $\low k\!\!\!k\!\!\!k^{(r)}_r$ is in fact $\low k^{(r)}_r$)
is the smallest
number in $T_i(\mu)$, thus it is the entry corresponding to the
first atypical root of $\mu$, \SSP i.e.,
\equa{MU-1}
{\mu_{\wt m_1}=\l_{m_r}-\low k\!\!\!k\!\!\!k^{(r)}_r.
\SSP }
We \SSP denote 
\equa{de-tau}
{\tau=(L^r_{m_r,n_r}(\l))^+=(\l-\low k\!\!\!k\!\!\!k^{(r)}_r d_{m_r,n_r})^+
\mbox{ \ \ (cf.~(\ref{low-nu}))}.
\SSP}
Then obviously we have
\VB
{\bf Claim 4}. The first atypical entry of $\tau$ is
$\mu_{\wt m_1}=\l_{m_r}-\low k\!\!\!k\!\!\!k^{(r)}_r$, and the $i$-th atypical entry of $\tau$
is $\l_{m_{i-1}}$ if $i>1$.
\VB
{\bf Claim 5}. $\wt k_1=\low k\!\!\!k\!\!\!k^{(r)}_r$
(this means that the element $\mu_{\wt m_1}+\theta'_1\wt k_1=\l_{m_r}$ is indeed
in the set $T(\l)$).
\VE
We prove this by induction on $\pi_1$ (cf.~(\ref{pi-i})).
If $\pi_1=1$, i.e., there does not exist $s<r$ with $a_s=1$.
Then conditions (\ref{conds4}) and (\ref{conds5}) imply
$c_{1,r}=\qq$, i.e., $\ell_{1,r}=r-1$, and $c_{s,r}=\nn$, i.e.,
$\ell_{s,r}>r-s$, for all $s>1.$
This together with (\ref{de-tau}) and (\ref{MU-1}) 
in turn 
implies that $\ell{\ssc\,}^\tau_{1,s}<s-1$ for all $s\in\II{2}{r}$
(recall Convention \ref{conv}).
Note 
from Claim 4 and (\ref{TTT}) that 
(recall Convention \ref{con5}\SSP) 
\equa{MAMANote}
{
\wt \ell_{1,s}\le\max\{\ell{\ssc\,}^\tau_{1,p}\,\,|\,\,p\in\II{1}{s}\}\le s-1.
\SSP}
This
shows that 
$\wt c_{1,s}=\cc$ for all $s>1$, i.e., $\wt p_1=r$. By (\ref{k-i-i}),
$\wt k_1=\wh{\wt k}_1^{(r)}$, which 
is indeed equal to $\low k\!\!\!k\!\!\!k^{(r)}_r$ from the definitions of $\low k\!\!\!k\!\!\!k^{(r)}_r$
and $\wh{\wt k}_1^{(r)}$ in (\ref{k---i}).

If $\pi_1>1$. We choose $p<r$ to be the largest such that $a_p=1$, then
$\mu_{\wt m_2}=\l_{m_p}-\low k\!\!\!k\!\!\!k^{(\theta_p)}_p$ (i.e., which is the second smallest
atypical entry of $\mu$) by the maximal choice of $p$, and one can prove that
$\wt c_{1,2}=\cc$ (i.e., $\wt\ell_{1,2}=0$). Using another induction on $\pi_1$,
one can again show that $\wt c_{1,s}=\cc$ for all $s>1$. Thus as above we have 
$\wt k_1=\low k\!\!\!k\!\!\!k^{(r)}_r$.
\VB
{\bf Claim 6}. If $p\le s<p+\pi_p$ for some $p$ (such $p$ must be unique), then
$\wt k_s=\low k\!\!\!k\!\!\!k^{(\theta_t)}_t$, where $t$ is the $(s+1-p)$-th
largest number such that $a_t=p$.
\VE
This claim follows from the proof of claim 5 
by forgetting all labels not in the set $\II{p}{t}$ 
(i.e., by restricting to $\II{p}{t}$)
and regarding
$\II{p}{t}$ as the whole set $\II{1}{r}$.
\VB
Now (\ref{suffice-prove}) follows from Claim 6, and thus the lemma.
\end{proof}
\begin{remark}\label{rem5.3-1}
The proof of Lemma \ref{lemm5.1} can also be done by induction on $|\theta|:=N=$ $
\#\{s\in\II{1}{r}\,|\,\theta_s\ne0\}$ in the following way:
The definition of $\tau$ in $(\ref{de-tau})$ shows that
$\mu=L'_{\theta^\tau}(\tau)$ for $\theta^\tau\in\Theta^\l$ with
$\theta^\tau_i=\theta_i$ if $i<r$ and $\theta^\tau_r=0$. So
by the inductive assumption, $\tau=R_{\theta''}(\mu)$ for some $\theta''\in\{0,1\}^r$
such that $\theta''_1=0$. Then one can obtain that $\l=R_{m_1,n_1}(\tau)
=R_{\theta'}(\mu)$, where
$\theta'_1=1$ and $\theta'_i=\theta''_i$ if $i>1$.

The arguments to be given below can be regarded as 
the ``reverse'' arguments of the above.
\end{remark}
Next suppose
$\l=R_{\theta'}(\mu)$ for some $\theta'\in\{0,1\}^r$.
First we note from the definition of (recall Convention \ref{con5})
$\wt p_s$ in (\ref{p-s}) \SSP that
\equa{p-s-}
{
s\le t\le \wt p_t\le \wt p_s\mbox{ \ for all \ }t\in\II{s}{\wt p_s}.
}
As in (\ref{subset}), this allows us to divide $\II{1}{r}=\VV_1\cup\VV_2$,
such that $\theta'_s=0$ if $s\in\VV_1$, and 
\equa{ell'-i-}
{
\theta'_{\ell'_i}=1\mbox{ \ and \ }
\wt p_{\ell'_1}=\ell_i,\,\ i=1,2,...,\nu.
}
(The above means that $\ell_i$ is the largest number
such that all atypical roots $\wt \g_s$, $\ell'_i<s\le\ell_i,$ of $\mu$ are $\cc$-related to
$\wt\g_{\ell'_i}$.)

As above we can suppose $\nu=1$ and $\ell'_i=1,\,\ell_i=r$, i.e.,
$\theta'_1=1,\,p_1=r$. 
This implies $\wt c_{1,s}=\cc$, \SSP i.e.,
\equa{c-1-s}
{
\wt \ell_{1,s}<s-1\mbox{ \ for all \ }s=2,3...,r.
\SSP}
We \SSP set
\equa{eta-}
{
\tau=R_{\wt m_1,\wt n_1}(\mu)=(\mu+\wt k_1d_{\wt m_1,\wt n_1})^+.
}
\begin{lemma}\label{ll-1}
The $i$-th atypical entry of $\tau$ is the $(i+1)$-th atypical entry
$\mu_{\wt m_{i+1}}$ if $i<r$, or 
$\mu_{\wt m_1,\wt n_1}+\wt k_1$ if $i=r$
$($cf.~Remarks \ref{addition}, \ref{rem5.3-1}$)$.
\end{lemma}
\begin{proof} This follows from (\ref{eta-}) and definition (\ref{k-i}).
\SSP\end{proof}
Set 
$\theta''\in\{0,1\}^r$ with 
$\theta''_i=\theta'_{i+1}$ if $i<r$ and $\theta''_r=0.$
This and definition (\ref{R'-theta-f}) \SSP imply 
\equa{SAA}
{
\l=R_{\theta'}(\mu)=R_{\theta''}(\tau).
}
By Lemma \ref{ll-1} and (\ref{eta-}), we obtain that for all $i<r$, $p^\tau_i<r$
(recall Convention \ref{conv}).
Since 
$|\theta''|:=\sum_{s=1}^r\theta''_s<|\theta'|$,
we can use induction on $|\theta'|$ to suppose
that
there exists some $\theta^\tau\in\Theta^\l$ with $\theta^\tau_r=0$, 
such that
\equa{SSSUCH}
{
\tau=L'_{\theta^\tau}(\l).
}
Note from 
definition (\ref{c-1-s})
that for all $s<r$ we must have $\ell{\ssc\,}^\tau_{s,r}\le r-1$, 
which together
with (\ref{SSSUCH}) in particular shows that
\equa{SAMMM}
{
\ell_{1,r}=\ell^\l_{1,r}\le \max\{\ell{\ssc\,}^\tau_{s,r}\,\,|\,\,s<r\}\le r-1.
}
Also (\ref{eta-}) (cf.~(\ref{SSSUCH})) implies that $\mu=L^r_{m_r,n_r}(\tau)$
(a lowering operator $L_{i,j}$ is the ``inverse operator'' of 
some raising operator $R_{i',j'}$ in some sense), which 
together with (\ref{SSSUCH}) and the definition of $L_\theta$ in
(\ref{l-theta-f}) \SSP implies
\equa{defined-by}
{\mu=L_\theta(\l),
\SSP}
where $\theta\in\N^r$ is defined by 
$\theta_i=\theta^\tau_i$ if $i<r$ and $\theta_r=r.$
\begin{lemma}\label{lemm-last}
$\theta\in\Theta^\l.$
\SSP\end{lemma}
\begin{proof}
The fact that $\theta^\tau\in\Theta^\l$ implies that conditions 
(\ref{conds1})--(\ref{conds5}) hold for all $s<r$. That these conditions
also hold for $s=r$
follows from the fact that
$\theta_r=r$ and $c_{1,r}\ne \nn$ (cf.~(\ref{SAMMM})).
\SSP\end{proof}
\ni{\it Proof of Theorem \ref{theo-main}}.
Finally we return to the proof of Theorem \ref{theo-main}. By
Lemma \ref{lemm5.1}, (\ref{defined-by}) and
Lemma \ref{lemm-last}, we have (\ref{to-prove}). 
This implies the theorem.
\hfill$\Box$

\vskip6pt
{\small\vskip4pt \noindent{\bf Acknowledgements}. The author
would like to thank Dr.~R.B.~Zhang for stimulated discussions.
This research was supported by
Australian Research Council, NSF grant 10171064 of
China, EYTP and TCTPFT grants of Ministry of Education of China.}



\begin{thebibliography}{9999}
\bibitem{BL} I.N.~Bernstein, D.A.~Leites,
{\em A formula for the characters of the irreducible finite-dimensional
representations of Lie superalgebras of series $gl$ and $sl$},
C.~R.~Acad.~Bulgare Sci.~{\bf33} (1980), 1049-1051.

\bibitem{B} J.~Brundan, {\em Kazhdan-Lusztig polynomials and
   character formulae for the Lie superalgebra $gl(m|n)$},
   J.~Amer.~Math.~Soc.~{\bf 16} (2002) 185-231.


\bibitem{CZ} S.-J.~Cheng, R.B.~Zhang, {\em  An analogue of
Kostant's ${\mathfrak u}$-cohomology
formula for the general linear superalgebra},
International Mathematics Research Notices {\bf1} (2004) 31-53.


\bibitem{Cr2} C.J.~Cummins, R.C.~King, {\em Composite Young diagrams, supercharacters of 
$U(M/N)$ and modification rules}, J.~Phys.~A {\bf20} (1987) 3121-3133.

\bibitem{G} J.~Germoni, {\em Indecomposable representations of special linear
Lie superalgebras}, J. Alg. {\bf209}, 367--401 (1998).

\bibitem{HKV} J.W.B.~Hughes, R.C.~King and J.~van der Jeugt,
{\em On the composition factors of Kac modules for the Lie superalgebras $sl(m|n)$},
J.~Math.~Phys.~{\bf33} (1992) 470-491.

\bibitem{Kac0} V.G.~Kac, {\em Classification of simple Lie superalgebras},
Funct.~Anal.~Appl. {\bf9} (1975) 263-265.

\bibitem{Kac} V.G.~Kac, {\em Lie superalgebras},
   Adv.~Math.~{\bf 26} (1977) 8-96.

\bibitem{Kac1} V.G.~Kac, {\em Characters of typical representations of 
  classical Lie superalgebras}, Comm.~Alg.~{\bf5} (1977) 889-897

\bibitem{Kac2} V.G.~Kac, {\em Representations of classical Lie superalgebras},
  Lect. Notes Math. {\bf676} (1978) 597-626.

\bibitem{P1} T.D.~Palev, {\em Finite-dimensional representations of the special linear
Lie superalgebra $sl(1,n)$ II. Nontypical representations},  
J.~Math.~Phys.~{\bf29} (1988) 2589-2598.

\bibitem{P2} T.D.~Palev, {\em Essentially typical representations of the Lie
 superalgebras $gl(n/m)$ in a Gel'fand-Zetlin basis},
Funct.~Anal.~Appl.~{\bf 23} (1989) 141-142 (English translation).



\bibitem{SZ98} M.~Scheunert, R.B.~Zhang,
   {\em Cohomology of Lie superalgebras and their generalizations},
   J.~Math.~Phys.~{\bf 39} (1998) 5024--5061.

\bibitem{SZ99} M.~Scheunert, R.B.~Zhang,
   {\em The second cohomology of $sl(m|1)$ with coefficients in its
    enveloping algebra is trivial},
   Lett.~Math.~Phys.~{\bf 47} (1999) 33--48.

\bibitem{SZ02} M.~Scheunert, R.B.~Zhang,
   {\em The general linear supergroup and its Hopf superalgebra
   of regular functions},
   J.~Algebra {\bf 254} (2002) 44--83.

\bibitem{Se96} V.~Serganova, {\em Kazhdan-Lusztig polynomials and
       character formula for the Lie superalgebra $\mathfrak{gl}(m|n)$},
       Selecta Math.~{\bf 2} (1996) 607-654.

\bibitem{Se98} V.~Serganova,
  {\em Characters of irreducible representations of simple Lie superalgebras},
  Proceedings of the International Congress of Mathematicians 1998, Berlin,
  Vol.~II, Documenta Mathematica, Journal der Deutschen
  Mathematiker-Vereinigung, pp.~583-593.



\bibitem{S} Y.~Su, {\em Weakly primitive vectors of
   Kac-modules of the Lie superalgebras $sl(m|n)$},
   J.~Math. Phys.~{\bf42} (2001) 5444-5456.

\bibitem{SHK} Y.~Su, J.W.B.~Hughes and R.C.~King, {\em Primitive vectors
  of Kac-modules of the Lie superalgebras $sl(m|n)$},
   J.~Math.~Phys.~{\bf41} (2000) 5064-5087.

\bibitem{SZ0} Y.~Su, R.B.~Zhang, {\em 
Cohomology of Lie superalgebras ${\mathfrak{sl}}_{m|n}$ and ${\mathfrak{osp}}_{2|2n}$},
preprint, math.QA/0402419.

\bibitem{SZ1} Y.~Su, R.B.~Zhang, {\em Character and dimension formulae for finite 
dimensional irreducible representations of the general linear superalgebra},
preprint, math.QA/0403315.

\bibitem{V} J.~van der Jeugt, {\em Character formulae for the 
  Lie superalgebra $C(n)$}, Comm.~Algebra {\bf19} (1991) 199-222.
 
\bibitem{VHKT0} J.~van der Jeugt, J.W.B.~Hughes, R.C.~King and J.~Thierry-Mieg,
{\em A character formula for singly atypical modules of the Lie superalgebra
$sl(m/n)$}, Comm.~Alg.~{\bf 19} (1991) 199-222.

\bibitem{VHKT} J.~van der Jeugt, J.W.B.~Hughes, R.C.~King and
J.~Thierry-Mieg, {\em Character formulas for irreducible modules of
the Lie superalgebras $sl(m|n)$}, J.~Math.~Phys.~{\bf31} (1990) 2278-2304.

\bibitem{VZ} J.~van der Jeugt, R.B.~Zhang, {\em Characters and composition
 factor multiplicities for the Lie superalgebra $gl(m|n)$},
  Lett.~Math.~Phys.~{\bf 47} (1999) 49-61.

\bibitem{Zh93} R.B.~Zhang,  {\em Finite dimensional irreducible representations of
the quantum supergroup} ${\rm U}_q({\mathfrak{gl}}(m/n))$,
J.~Math.~Phys.~{\bf 34} (1993) 1236--1254.

\bibitem{Zh98} R.B.~Zhang,
   {\em Structure and representations of the quantum general linear
   supergroup}, Commun. Math.~Phys.~{\bf 195} (1998) 525--547.

\bibitem{Zou} Y.M.~Zou, {\em Categories of finite-dimensional weight
modules over type I classical Lie superalgebras},  J.~Algebra {\bf 180} (1996) 459-482.

\end{thebibliography}
\end{document}